\nonstopmode
\documentclass[]{amsart}
\usepackage{amssymb}
\usepackage{amscd}
\usepackage{graphicx}
\usepackage{psfrag}
\usepackage{epsfig}
\usepackage{verbatim}% for comment enviroment
\def\undersetbrace#1\to#2{\underbrace{#2}_{#1}}
\def\oversetbrace#1\to#2{\overbrace{#2}^{#1}}
\def\AMSunderset#1\to#2{\underset{#1}{#2}}
\def\AMSoverset#1\to#2{\overset{#1}{#2}}
\def\therosteritem#1{{\rm (#1)}}
\def\East#1#2{\overset{#1}{\longrightarrow}}
\newenvironment{proclaim}[1]{\par\medskip\noindent{\bf #1.}\it}{\par\smallskip}
\newenvironment{demo}[1]{\par\smallskip\noindent{\bf #1.}}{\par\smallskip}
\def\nmb#1#2{\if0#1{#2}%
\else\if:#1{#2}%
\else\if|#1{{#2}}
\else\if.#1{\gdef\lastnmb{#2}#2}%
\else{{#2}}%
\fi\fi\fi\fi}
\def\AMStag#1{\tag{#1}}
\def\cit#1#2{\ifx#1!\cite{#2}\else#2\fi} %for citing references
\def\totoc{}             %= to table of content, invoked by kms-book.sty
\def\idx{}               % for producing index, invoked by kms-book.sty
\def\ign#1{}             %=ignore, invisible entry for the index only
\def\o{\circ}

\def\al{\alpha}
\def\be{\beta}
\def\ga{\gamma}
\def\de{\delta}
\def\ep{\varepsilon}

\def\th{\theta}

\def\ka{\kappa}

\def\si{\sigma}

\def\ph{\varphi}

\def\ps{\psi}
\def\om{\omega}
\def\Ga{\Gamma}

\def\Ph{\Phi}

\def\Om{\Omega}
\def\i{^{-1}}
\def\x{\times}
\def\p{\partial}
\let\on=\operatorname
\def\today{\ifcase\month\or
 January\or February\or March\or April\or May\or June\or
 July\or August\or September\or October\or November\or December\fi
 \space\number\day, \number\year}
\def\AMSonly#1{}

\begin{document}%\topmatter

\title{Riemannian geometries on spaces of plane curves
}
\author{Peter W. Michor, David Mumford}
\address{
Peter W. Michor:
Fakult\"at f\"ur Mathematik, Universit\"at Wien,
Nordberg\-strasse 15, A-1090 Wien, Austria; {\it and:}
Erwin Schr\"odinger Institut f\"ur Mathematische Physik,
Boltzmanngasse 9, A-1090 Wien, Austria
}
\email{Peter.Michor@univie.ac.at}
\address{
David Mumford:
Division of Applied Mathematics, Brown University,
Box F, Providence, RI 02912, USA}
\email{David\_{}Mumford@brown.edu}
%\begin{comment}\dedicatory \enddedicatory\end{comment}
\date{\today}
\thanks{Both authors were supported by NSF Grant 007~4276.
     P.W.M. was supported by FWF Project P~14195~MAT'}

%\keywords{}
\begin{abstract} We study some Riemannian metrics on the space of smooth regular 
curves in the plane, viewed as the orbit space of maps from $S^1$ 
to the plane modulo the group of diffeomorphisms of $S^1$, acting as 
reparameterizations. In particular we investigate the metric
for a constant $A> 0$:
\begin{displaymath}
G^A_c(h,k) := \int_{S^1}(1+A\ka_c(\th)^2)\langle h(\th),k(\th) \rangle 
     |c'(\th)|\,d\th
\end{displaymath}
where $\ka_c$ is the curvature of the curve $c$ and $h,k$ are normal vector
fields to $c$. The term $A\ka^2$ is a sort of geometric Tikhonov regularization 
because, for $A=0$, the geodesic distance between any 2 distinct curves is 0, 
while for $A>0$ the distance is always positive. We give some lower bounds
for the distance function, derive the geodesic equation and the sectional
curvature, solve the geodesic equation with simple endpoints numerically,
and pose some open questions. The space has an interesting split personality: 
among large smooth curves, all its sectional curvatures are $\ge 0$,
while for curves with high curvature or perturbations of high frequency,
the curvatures are $\le 0$.
\end{abstract}
\def\LaTeXonly{}%\endtopmatter
%%\input amspptb.sty
%\userunningheads
%\def%\leftheadtext{\scshape Peter W. Michor}
%\def\rightheadtext{\scshape }
%\def\bottremark{\today\hfill}

\maketitle

%\section*{Table of contents }
%%\input \jobname.toc
%\tableofcontents
%\makeindex

\section*{\totoc\nmb0{1}. Introduction }

This paper arose from the attempt to find the simplest Riemannian metric on the
space of 2-dimensional `shapes'. By a shape we mean a compact simply connected 
region in the plane whose boundary is a simple closed curve. By requiring that
the boundary curve has various degrees of smoothness, we get not just one space
but a whole hierarchy of spaces. All these spaces will include, however, a
core, namely the space of all shapes with $C^\infty$ boundary curves. We expect 
that the most natural shape spaces will arise as the completions of this core 
space in some metric hence we take this core as our basic space. Note that it is
the orbit space 
\begin{displaymath}
B_e(S^1,\mathbb R^2)=\on{Emb}(S^1,\mathbb R^2)/\on{Diff}(S^1)\end{displaymath}
of the space of all $C^\infty$ embeddings of $S^1$ in the plane, under the
action by composition from the right by diffeomorphisms of the circle. 
The space $\on{Emb}(S^1,\mathbb R^2)$ is a smooth manifold, in fact an open
subset of the Fr\'echet space $C^\infty(S^1,\mathbb R^2)$, and it is the total 
space of a smooth principal bundle with base $B_e(S^1,\mathbb R^2)$

In fact, most of our results carry over to the bigger orbit space of 
immersions mod diffeomorphisms:
\begin{displaymath}
B_i(S^1, \mathbb R^2) =\on{Imm}(S^1, \mathbb R^2)/\on{Diff}(S^1).\end{displaymath}
This action is not quite free (see \nmb!{2.4} and \nmb!{2.5}), 
hence this orbit space is an 
orbifold (see \nmb!{2.5}) and not quite a manifold. 
There is the slightly smaller space $\on{Imm}_{f}$ (see \nmb!{2.1})
of immersed curves where diffeomorphisms act freely, the total space of a
principal fiber bundle with a natural connection admitting parallel
transport. Existence of horizontal curves, however, holds also in the
big space $\on{Imm}$ (see \nmb!{2.5}) which will be one of the weapons in
our hunt for geodesics on $B_i$.

The second author was led to study the space $B_e$ from its relevance to 
computer vision. To understand an image of the world, one needs to identify the 
most salient objects present in this image. In addition to readily quantifiable
properties like color and area, objects in the world and their projections 
depicted by 2D images possess a `shape' which is readily used by human
observers to distinguish, for example, cats from dogs, BMW's from Hondas, etc.
In fact people are not puzzled by what it means to say two shapes are 
{\it similar} but rather find this a natural question. This suggests that we 
construct, on some crude level, a mental metric which can be used to recognize 
familiar objects by the similarity of their shapes and to cluster categories of 
related objects like cats. Incidentally, immersions also arise in 
vision when a 3D object partially occludes itself from some viewpoint, hence 
its full 2D contour has visible and invisible parts which, together, form 
an immersed curve in the image plane.

It is a central problem in computer vision to devise algorithms by which
computers can similarly recognize and cluster shapes. Many types of metrics 
have been proposed for this purpose \cit!{7}. For example, 
there are $L^1$-type metrics such are the area of the symmetric difference of the 
interiors of two shapes. And there are $L^\infty$-type
metrics such as the Hausdorff metric: the maximum distance of points on 
either shape from the points on the other or of points outside one shape from
points outside the other. These metrics will come up below,
but the starting point of this investigation was whether one could use the
manifold structure on the space of shapes and define an $L^2$-type metric by
introducing a Riemannian structure on the space. 

Such questions have also arisen in Teichm\"uller theory and string theory, where
the so-called Weil-Peterssen metric on the space of shapes (also called 
the `universal Teich\-m\"uller space) has been much studied. In a second part of
this paper, we will compare our metric to this remarkable (homogeneous!) metric.

In this paper, we sought the absolutely simplest Riemannian metric that the
space $B_i$ supports. The most obvious $\on{Diff}(S^1)$-invariant weak Riemannian 
metric on the space of immersions is the $H^0$-metric:
\begin{displaymath}
G^0_c(h,k)=\int_{S^1}\langle h(\th),k(\th) \rangle |c'(\th)| \,d\th
\end{displaymath}
where $c:S^1 \rightarrow \mathbb R^2$ is an embedding defining a point in $B_e$ 
and $h,k$ are vector fields along the image curve, defining two tangent
vectors to $\on{Imm}(S^1,\mathbb R^2)$ at $c$. 
This induces a $\on{Diff}(S^1)$-invariant weak Riemannian metric on the 
space of all immersions and on $\on{Emb}(S^1,\mathbb R^2)$, 
and for the latter space it
induces a weak Riemannian metric on the base manifold $B_e$.

Surprisingly, the Riemannian distance defined as the infimum of the 
arclength of paths connecting two points in $B_e(S^1,\mathbb R^2)$ turns out to 
be 0, see \nmb!{3.10}! This seems to be one of the first examples where this
purely infinite dimensional phenomenon actually appears. 

Motivated by the proof of this result \nmb!{3.10} we are led to consider the
invariant Riemannian metric \nmb!{3.2.6} for a constant $A> 0$:
\begin{displaymath}
G^A_c(h,k) := \int_{S^1}(1+A\ka_c(\th)^2)\langle h(\th),k(\th) \rangle 
     |c'(\th)|\,d\th
\end{displaymath}
where $\ka_c(\th)$ is the curvature of $c$ at $c(\th)$.
We will argue that this induces a reasonable metric on $B_e(S^1,\mathbb R^2)$,
as the infimum of arclengths of paths connecting distinct points is always 
positive. Another reason is that the length function 
$\ell:B_e(S^1,\mathbb R^2)\to\mathbb R_{\ge0}$ has the following Lipschitz 
estimate \nmb!{3.3.2} with respect to this Riemannian distance:
\begin{displaymath}
\sqrt{\ell(C_1)}-\sqrt{\ell(C_0)}\le
\frac1{2\sqrt{A}}\on{dist}^{B_e}_{G^A}(C_1,C_2).
\end{displaymath}
In fact, one can bound the Fr\'echet distance between two curves in terms of this
metric (see \nmb!{3.5}). The completion of the space of smooth curves in this
metric contains all curves for whose curvature exists weakly as a finite
signed measure (e.g.\ piecewise $C^2$ curves) and is contained in the space 
of Lipschitz maps from $S^1$ to $\mathbb R^2$ modulo a suitable equivalence 
relation, see \nmb!{2.11}.

The geodesic equation for the metric $G^A$ on $\on{Emb}(S^1,\mathbb R^2)$ 
and on $B_e(S^1,\mathbb R^2)$
can be found in \nmb!{4.1.1}: It is a highly non-linear partial
differential equation of order 4 with degenerate symbol, but which nonetheless
seems to have a hypoelliptic linearization. 
If $A=0$, the equation reduces to a non-linear second order hyperbolic PDE,
which gives a well defined local geodesic spray. For any $A$, the sectional 
curvature on $B_e(S^1,\mathbb R^2)$ has an elegant expression which 
can be found in \nmb!{4.6.2} and \nmb!{4.6.4}. It is non-negative if $A=0$
and, for general $A$, becomes strictly negative only if the curve has large
curvature or the plane section has high frequency.
Of course we would have liked to solve the problem of existence and
uniqueness of geodesics for $A>0$. We can, however, translate the minimization 
of path length in our metric into an anisotropic Plateau-like problem:
In \nmb!{3.12} we show that a curve projects onto a geodesic in
$B_e(S^1,\mathbb R^2)$ if and only if its graph in $[0,1]\x \mathbb R^2$ is a
surface with given boundary at $\{0\}\x \mathbb R^2$ and $\{1\}\x \mathbb R^2$
which is critical for the anisotropic area functional \nmb!{3.12.3}.

In \nmb!{5.1} we determine the geodesic running through concentric circles
and the equation for Jacobi vector fields along this geodesic. The solution of 
the ordinary differential equation \nmb!{5.1.1} describing this geodesic can be 
written in terms of elliptic functions. This geodesic is no longer globally 
minimizing when the radius of the circles is large compared
to $\sqrt{A}$ and has conjugate points when it hits this positive curvature 
zone. 
In \nmb!{5.2} we study geodesics connecting arbitrary distant curves, 
hence requiring long translations. The middle part of such geodesics appear to 
be approximated by a uniformly translating `cigar'-like curve with semi-circular 
ends of radius $\sqrt{A}$ connected by straight line segments parallel to the
direction of translation. These figures were found by numerically minimizing a 
discrete form of the energy functional \nmb!{3.12.1}.

Finally, in \nmb!{5.3} and \nmb!{5.4}, we have some further pictures of 
geodesics. First we examine the formation of singularities when a small
perturbation is propagated forward and $A=0$. Then we look at some
geodesic triangles in $B_e$ whose vertices are ellipses with the same
eccentricity and center but different orientations. For various values
of $A$, we find that these triangles have angle sums greater and less than $\pi$.

\section*{\totoc\nmb0{2}. The manifold of immersed closed curves
}

\subsection*{\nmb.{2.1}. Conventions }

It is often convenient to use the identification $\mathbb R^2\cong\mathbb C$, 
giving us:
\begin{displaymath}
\bar x y = \langle x,y \rangle + i\det(x,y),\quad
\det(x,y)=\langle ix,y\rangle.
\end{displaymath}
We shall use the following spaces of $C^\infty$ (smooth) 
diffeomorphisms and curves, and we give the shorthand and the full name:

\begingroup
\parindent=0cm
\everypar{\hangindent=1.5cm \hangafter=1} 

$\on{Diff}(S^1)$, the regular Lie group (\cit!{4},~38.4) of all
diffeomorphisms $S^1\to S^1$ with its connected components
$\on{Diff}^+(S^1)$ of orientation preserving diffeomorphisms and 
$\on{Diff}^-(S^1)$ of orientation reversing diffeomorphisms. 

$\on{Diff}_1(S^1)$, the subgroup of diffeomorphisms fixing $1\in S^1$.
We have diffeomorphically 
$\on{Diff}(S^1)=\on{Diff}_1(S^1)\x S^1=\on{Diff}_1^+(S^1)\x
(S^1\rtimes \mathbb Z_2)$.

$\on{Emb}=\on{Emb}(S^1,\mathbb R^2)$, 
the manifold of all smooth embeddings $S^1\to \mathbb R^2$. 
Its tangent bundle is given by 
$T\on{Emb}(S^1,\mathbb R^2)=\on{Emb}(S^1,\mathbb R^2)\x C^\infty(S^1,\mathbb R^2)$.

$\on{Imm}=\on{Imm}(S^1,\mathbb R^2)$, 
the manifold of all smooth immersions $S^1\to \mathbb R^2$.
Its tangent bundle is given by 
$T\on{Imm}(S^1,\mathbb R^2)=\on{Imm}(S^1,\mathbb R^2)\x C^\infty(S^1,\mathbb R^2)$.

$\on{Imm}_{f}=\on{Imm}_{f}(S^1,\mathbb R^2)$, 
the manifold of all smooth 
free immersions $S^1\to \mathbb R^2$, i.e., those with trivial isotropy group
for the right action of $\on{Diff}(S^1)$ on $\on{Imm}(S^1,\mathbb R^2)$.

$B_e=B_e(S^1,\mathbb R^2)=\on{Emb}(S^1,\mathbb R^2)/\on{Diff}(S^1)$, the manifold of
1-dimensional connected submanifolds of $\mathbb R^2$, see \nmb!{2.3}.

$B_i=B_{i}(S^1,\mathbb R^2)=\on{Imm}(S^1,\mathbb R^2)/\on{Diff}(S^1)$, an infinite
dimensional `orbifold'; its points are, roughly speaking, smooth curves
with crossings and multiplicities, see \nmb!{2.5}.

$B_{i,f}=B_{i,f}(S^1,\mathbb R^2)=\on{Imm}_{f}(S^1,\mathbb
R^2)/\on{Diff}(S^1)$, a manifold, the base of a principal fiber bundle, see
\nmb!{2.4.3}.

\endgroup

We want to avoid referring to a path in our infinite dimensional spaces like 
$\on{Imm}$ or $B_e$ as a curve, 
because it is then a `curve of curves' and confusion arises
when you refer to a curve. So we will always talk of {\it paths} in the
infinite dimensional spaces, not curves. Curves will be in $\mathbb R^2$.
Moreover, if $t\mapsto (\th\mapsto c(t,\th))$ is a path, its $t$-th curve
will be denoted by $c(t)=c(t,\quad)$. By $c_t$ we shall denote the
derivative $\p_t c$, and $c_\th =\p_\th c$.  

\subsection*{\nmb.{2.2}. Length and curvature on 
$\on{Imm}(S^1,\mathbb R^2)$ } 

The volume form on $S^1$ induced by $c$ is given by 
\begin{equation}
\on{vol}:\on{Emb}(S^1,\mathbb R^2) \to \Om^1(S^1),\qquad \on{vol}(c) = 
|c_\th|\,d\th
\AMStag{\nmb:{1}}\end{equation}
and its derivative is 
\begin{equation}
d\on{vol}(c)(h) = \frac{\langle h_\th,c_\th \rangle}{|c_\th|}d\th.
\AMStag{\nmb:{2}}\end{equation}
We shall also use the {\it normal unit field}
\begin{displaymath}
n_c = i\frac{c_\th}{|c_\th|}.
\end{displaymath}
The length function is given by 
\begin{equation}
\ell:\on{Imm}(S^1,\mathbb R^2) \to \mathbb R,\qquad \ell(c) 
  = \int_{S^1} |c_\th|\,d\th
\AMStag{\nmb:{3}}\end{equation}
and its differential is  
\begin{align}
d\ell(c)(h) &= \int_{S^1} \frac{\langle h_\th,c_\th \rangle}{|c_\th|}d\th 
=-\int_{S^1} \Bigl\langle h, \frac{c_{\th\th}}{|c_\th|} 
-\frac{\langle c_{\th\th},c_\th\rangle}{|c_\th|^3}c_\th \Bigr\rangle d\th
\AMStag{\nmb:{4}}\\&
=-\int_{S^1} \langle h, \ka(c).ic_\th\rangle d\th
=-\int_{S^1} \langle h, n_c\rangle\ka(c) \on{vol}(c)
\notag
\end{align}

The curvature mapping is given by 
\begin{equation}
\ka:\on{Imm}(S^1,\mathbb R^2)\to C^\infty(S^1,\mathbb R),\qquad
\ka(c)= \frac{\det(c_\th,c_{\th\th})}{|c_\th|^3}
= \frac{\langle  ic_\th,c_{\th\th}\rangle}{|c_\th|^3}
\AMStag{\nmb:{5}}\end{equation}
and is equivariant so that $\ka(c\o f)=\pm\ka(c)\o f$ for 
$f\in \on{Diff}^\pm(S^1)$.
Its derivative is given by 
\begin{equation}
d\ka(c)(h) = \frac{\langle ih_\th,c_{\th\th} \rangle}{|c_\th|^3}
     + \frac{\langle ic_\th,h_{\th\th} \rangle}{|c_\th|^3}
     -3\ka(c)\frac{\langle h_\th,c_\th \rangle}{|c_\th|^2}.
\AMStag{\nmb:{6}}\end{equation}
With some work, this can be shown to equal:
\begin{equation}\AMStag{\nmb:{7}} 
d\ka(c)(h) = \frac{\langle h, c_\th \rangle}{|c_\th|^2} \ka_\th
+ \frac{\langle h, ic_\th \rangle}{|c_\th|}\ka^2
+ \frac{1}{|c_\th|}\Bigl( \frac{1}{|c_\th|} 
\bigl( \frac{\langle h, ic_\th \rangle}{|c_\th|}  \bigr)_\th \Bigr)_\th.
\end{equation}
To verify this, note that both the left and right hand side are equivariant with
respect to $\on{Diff}(S^1)$, hence it suffices to check it for constant
speed parametrizations, i.e.\ $|c_\th|$ is constant and
$c_{\th \th} = \ka |c_\th| i c_\th$. By linearity, it is enough
to take the 2 cases $h = aic_\th$ and $h=bc_\th$. Substituting 
these into formulas \thetag{\nmb|{6}} and \thetag{\nmb|{7}}, 
the result is straightforward.

\subsection*{\nmb.{2.3}. The principal bundle of embeddings 
$\on{Emb}(S^1,\mathbb R^2)$}

We recall some basic results whose proof can be found in {\cit!{4}:}

{\it \noindent (A) The set $\on{Emb}(S^1,\mathbb R^2)$ of all smooth embeddings 
$S^1\to \mathbb R^2$ is an open subset of the Fr\'echet space 
$C^\infty(S^1,\mathbb R^2)$ of all smooth mappings $S^1\to \mathbb R^2$ 
with the $C^\infty$-topology. 
It is the total space of a smooth principal bundle
$\pi:\on{Emb}(S^1,\mathbb R^2)\to B_e(S^1,\mathbb R^2)$ with structure group 
$\on{Diff}(S^1)$, the smooth regular Lie group group of all 
diffeomorphisms of $S^1$, whose base $B_e(S^1,\mathbb R^2)$ is the smooth 
Fr\'echet manifold of all submanifolds of $\mathbb R^2$ of type $S^1$, 
i.e., the smooth manifold of all simple closed curves in $\mathbb R^2$.
{\rm (\cit!{4},~44.1)} 

\noindent (B) This principal bundle admits a smooth principal connection
described by the horizontal bundle whose fiber $\mathcal N_c$ over $c$ 
consists of all vector fields $h$ along $c$ such that
$\langle h,c_\th \rangle=0$. The parallel transport for this connection
exists and is smooth. {\rm (\cit!{4},~39.1 and 43.1)}}

See \nmb!{2.4.3} for a sketch of proof of the first part in a slightly more
general situation. 
See also \nmb!{3.2.2} and \nmb!{3.2.3} for the horizontal bundle 
$\mathcal N_c$.
Here we want to sketch the use of the second part.
Suppose that $t\mapsto (\th\mapsto c(t,\th))$ is a path in 
$\on{Emb}(S^1,\mathbb R^2)$. Then $\pi\o c$ is a smooth path in $B_e(S^1,\mathbb
R^2)$. Parallel transport over it with initial value $c(0,\cdot)$ is a now
a path $f$ in $\on{Emb}(S^1,\mathbb R^2)$ which is horizontal, i.e., we have
$\langle f_t,f_\th \rangle=0$. This argument will play an important role
below. In \nmb!{2.5} below we will prove this property for general
immersions.

\subsection*{\nmb.{2.4}. Free immersions
}

The manifold $\on{Imm}(S^1,\mathbb R^2)$ of all immersions $S^1\to \mathbb R^2$
is an open set in the manifold $C^\infty(S^1,\mathbb R^2)$ and thus itself a
smooth manifold. 
An immersion $c:S^1\to \mathbb R^2$ is called {\it free} if $\on{Diff}(S^1)$
acts freely on it, i.e., $c\o \ph=c$ for $\ph\in \on{Diff}(S^1)$ implies
$\ph=\on{Id}$. 
We have the following results:

\noindent\therosteritem{\nmb:{1}}
{\it If $\ph\in \on{Diff}(S^1)$ has a fixed point and if $c\o\ph=c$ for
some immersion $c$ then $\ph=\on{Id}$.} This is (\cit!{1},~1.3).

\noindent\therosteritem{\nmb:{2}}
{\it If for $c\in \on{Imm}(S^1,\mathbb R^2)$ there is a point $x\in c(S^1)$
with only one preimage then $c$ is a free immersion.} This is
(\cit!{1},~1.4). There exist free immersions without such points: Consider
a figure eight consisting of two touching ovals, and map $S^1$ to this by
first transversing the upper oval 3 times and then the lower oval 2
times. This is a free immersion. 

\begin{proclaim}{\therosteritem{\nmb:{3}} The manifold $B_{i,f}(S^1,\mathbb R^2)$}
{\rm (\cit!{1},~1.5)}
The set $\on{Imm}_{f}(S^1,\mathbb R^2)$ of all free immersions
is open in $C^\infty(S^1,\mathbb R^2)$ and thus a smooth submanifold. 
The projection 
\begin{displaymath}
\pi:\on{Imm}_{f}(S^1,\mathbb R^2)\to 
\frac{\on{Imm}_{f}(S^1,\mathbb R^2)}{\on{Diff}(S^1)}=:B_{i,f}(S^1,\mathbb R^2)
\end{displaymath}
onto a Hausdorff smooth manifold is a smooth principal fibration with structure
group $\on{Diff}(S^1)$. 
By {\rm (\cit!{4},~39.1 and 43.1)} this fibration 
admits a smooth principal connection described by the
horizontal bundle with fiber $\mathcal N_c$ consisting of all vector fields $h$ 
along $c$ such that $\langle h,c_\th \rangle=0$. This connection admits a
smooth parallel transport over each each smooth curve in the base manifold.
\end{proclaim}

We might view $\on{Imm}_{f}(S^1,\mathbb R^2)$ as the nonlinear Stiefel 
manifold of parametrized curves in $\mathbb R^2$ and consequently 
$B_{i,f}(S^1,\mathbb R^2)$ as the nonlinear Grassmannian of unparametrized 
simple closed curves. 

\begin{demo}{Sketch of proof} 
See also \cit!{1} for a slightly different
proof with more details.
For $c\in \on{Imm}_{f}(S^1,\mathbb R^2)$ and
$s=(s_1,s_2)\in \mathcal V(c)\subset C^\infty(S^1,\mathbb R\x S^1)$ consider 
\begin{displaymath}
\ph_c(s):S^1\to \mathbb R^2,\qquad \ph_c(s)(\th)=c(s_2(\th))+s_1(s_2(\th)).n_c(s_2(\th))
\end{displaymath}
where $\mathcal V(c)$ is a $C^\infty$-open neighborhood of $(0,\on{Id}_{S^1})$
in $C^\infty(S^1,\mathbb R\x S^1)$
chosen in such a way that: 
\begin{enumerate}
\item[$\bullet$] 
  $s_2\in\on{Diff}(S^1)$ for each $s\in\mathcal V(c)$.
\item[$\bullet$] 
  $\ph_c(s)$ is a free immersion for each $s\in\mathcal V(c)$.
\item[$\bullet$] 
  For $(s_1,s_2)\in \mathcal V(c)$ and $\al\in\on{Diff}(S^1)$ we have
  $(s_1,s_2\o\al)\in \mathcal V(c)$.
\end{enumerate}
Obviously $\ph_c(s_1,s_2)\o\al=\ph_c(s_1,s_2\o\al)$ and $s_2$ is uniquely
determined by $\ph_c(s_1,s_2)$ since this is a free immersion.
Thus the inverse of $\ph_c$ is a
smooth chart for the manifold $\on{Imm}_{f}(S^1,\mathbb R^2)$.
Moreover, we consider the mapping (which will be important in section \nmb!{4}
below)
\begin{align*}
&\ps_c: C^\infty(S^1,(-\ep,\ep)) \to \on{Imm}_{f}(S^1,\mathbb R^2),
  \qquad \mathcal Q(c):= \ps_c(C^\infty(S^1,(-\ep,\ep)))\\
&\ps_c(f)(\th) = c(\th) + f(\th)n_c(\th) = \ph_c(f,\on{Id}_{S^1})(\th),\\
&\pi\o\ps :C^\infty(S^1,(-\ep,\ep)) \to B_{i,f}(S^1,\mathbb R^2),
\end{align*}where $\ep$ is small. Then (an open subset of) $\mathcal V(c)$ splits
diffeomorphically into 
\begin{displaymath}
C^\infty(S^1,(-\ep,\ep)) \x \on{Diff}{S^1}\end{displaymath}
and thus its image under $\ph_c$ splits into $\mathcal Q(c)\x \on{Diff}(S^1)$.
So the inverse of $\pi\o\ps_c$ is a smooth chart for 
$B_{i,f}(S^1,\mathbb R^2)$.
That the chart changes induced by the mappings $\ph_c$ and $\ps_c$ 
contructed here are smooth is shown 
by writing them in terms of compositions and projections only 
and applying the setting of \cit!{4}.
\qed\end{demo}

\subsection*{\nmb.{2.5}. Non free immersions
}

Any immersion is proper since $S^1$ is compact and thus 
by (\cit!{1},~2.1) the orbit space 
$B_i(S^1,\mathbb R^2)=\on{Imm}(S^1,\mathbb R^2)/\on{Diff}(S^1)$ is Hausdorff. 
Moreover, by (\cit!{1},~3.1 and 3.2) for any immersion $c$ 
the isotropy group $\on{Diff}(S^1)_c$ 
is a finite cyclic group which acts as group of covering transformations for a
finite covering $q_c:S^1\to S^1$ such that $c$ factors over $q_c$ to a free
immersion $\bar c:S^1\to \mathbb R^2$ with $\bar c\o q_c = c$. 
Thus the subgroup $\on{Diff}_1(S^1)$ of all diffeomorphisms $\ph$ fixing
$1\in S^1$ acts freely on $\on{Imm}(S^1,\mathbb R^2)$.
Moreover, for each $c\in \on{Imm}$  the submanifold
$\mathcal Q(c)$ from the proof of \nmb!{2.4.3} (dropping the freeness
assumption) is a slice in a strong sense:
\begin{enumerate}
\item[$\bullet$] 
   $\mathcal Q(c)$ is invariant under the isotropy group $\on{Diff}(S^1)_c$.
\item[$\bullet$] 
   If $\mathcal Q(c)\o \ph \cap \mathcal Q(c)\ne\emptyset$ for
   $\ph\in\on{Diff}(S^1)$ then $\ph$ is already in the isotropy group 
   $\ph\in\on{Diff}(S^1)_c$.
\item[$\bullet$] 
   $\mathcal Q(c)\o \on{Diff}(S^1)$ is an invariant open neigbourhood of
   the orbit $c\o\on{Diff}(S^1)$ in $\on{Imm}(S^1,\mathbb R^2)$ which
   admits a smooth retraction $r$ onto the orbit. The fiber $r\i(c\o\ph)$
   equals $\mathcal Q(c\o\ph)$.
\end{enumerate}
Note that also the action 
$$
\on{Imm}(S^1,\mathbb R^2)\x \on{Diff}(S^1) \to 
\on{Imm}(S^1,\mathbb R^2)\x\on{Imm}(S^1,\mathbb R^2),\qquad (c,\ph)\mapsto
(c,c\o\ph)
$$
is proper so that all assumptions and conclusions of Palais' slice theorem
\cit!{P} hold.
This results show that the orbit space 
$B_i(S^1,\mathbb R^2)$ has only very simple singularities of the type of a cone
$\mathbb C/\{e^{2\pi k/n}: 0\le k<n\}$ times a Fr\'echet space.
We may call the space $B_i(S^1,\mathbb R^2)$ an infinite dimensional {\it
orbifold}. The projection 
$\pi:\on{Imm}(S^1,\mathbb R^2)\to B_i(S^1,\mathbb R^2)
=\on{Imm}(S^1,\mathbb R^2)/\on{Diff}(S^1)$
is a submersion off the singular points and has only mild singularities at the
singular strata.
The normal bundle $\mathcal N_c$ mentioned in \nmb!{2.3} is well defined
and is a smooth vector subbundle of the tangent bundle. 
We do not have a principal bundle and thus no principal connections, but we
can prove the main consequence, the existence of horizontal paths,
directly:

\begin{proclaim}{Proposition}
For any smooth path $c$ in $\on{Imm}(S^1,\mathbb R^2)$ there exists a smooth
path $\ph$ in $\on{Diff}(S^1)$ with $\ph(0,\quad)=\on{Id}_{S^1}$ depending
smoothly on $c$ such that
the path $e$ given by $e(t,\th)=c(t,\ph(t,\th))$ is horizontal:
$e_t\bot e_\th$.
\end{proclaim}

\begin{demo}{Proof}
Let us write $e=c\o \ph$ for $e(t,\th)=c(t,\ph(t,\th))$, etc. 
We look for $\ph$
as the integral curve of a time dependent vector field $\xi(t,\th)$ on
$S^1$, given by $\ph_t=\xi\o \ph$.
We want the following expression to vanish:
\begin{align*}
\langle \p_t(c\o\ph),\p_\th(c\o\ph) \rangle 
&=\langle c_t\o\ph + (c_\th\o\ph)\,\ph_t,(c_\th\o\ph)\,\ph_\th \rangle 
\\&
=(\langle c_t,c_\th\rangle\o\ph)\,\ph_\th +(\langle c_\th,
c_\th\rangle\o\ph)\,\ph_\th\,\ph_t
\\
&=\bigl((\langle c_t,c_\th\rangle +\langle c_\th,
c_\th\rangle\,\xi)\o\ph\bigr)\,\ph_\th.
\end{align*}
Using the time dependent vector field $\xi=-\frac{\langle
c_t,c_\th\rangle}{|c_\th|^2}$ and its flow $\ph$ achieves this.
\qed\end{demo}

\subsection*{\nmb.{2.6}. The manifold of immersions with constant speed}
Let $\on{Imm}_{a}(S^1,\mathbb R^2)$ be the space of all immersions
$c:S^1\to \mathbb R^2$ which are parametrized by scaled arc length, so that
$|c_\th|$ is constant. 

\begin{proclaim}{Proposition} 
The space $\on{Imm}_{a}(S^1,\mathbb R^2)$ is a smooth manifold. There is a
diffeomorphism 
$\on{Imm}(S^1,\mathbb R^2)=\on{Imm}_{a}(S^1,\mathbb R^2)\x \on{Diff}_1^+(S^1)$
which respects the splitting $\on{Diff}(S^1)=\on{Diff}_1^+(S^1)\ltimes
(S^1\ltimes \mathbb Z_2)$. 
There is a smooth action of the rotation and reflection
group $S^1\ltimes \mathbb Z_2$ on
$\on{Imm}_{a}(S^1,\mathbb R^2)$ with orbit space
$\on{Imm}_{a}(S^1,\mathbb R^2)/(S^1\ltimes\mathbb Z_2)=B_i(S^1,\mathbb R^2)$.  
\end{proclaim}

\begin{demo}{Proof}
For $c\in \on{Imm}(S^1,\mathbb R^2)$ we 
put
\begin{align*}
&\si_{c}\in\on{Diff}_1(S^1),\qquad\qquad
\si_{c}(\th) = \exp\Bigl(\frac
     {2\pi i\int_1^{\th}|c'(u)|du}{\int_{S^1}|c'(u)|du}\Bigr) 
\\
&\al:\on{Imm}(S^1,\mathbb R^2)\to \on{Imm}_a(S^1,\mathbb R^2),\qquad  
\al(c)(\th) := c(\si_{c}\i(\th)). 
\end{align*}
By the fundamentals of manifolds of mappings \cit!{4} the mapping 
$\al$ is smooth from $\on{Imm}(S^1,\mathbb R^2)$ into itself and we have 
$\al\o\al=\al$. 

Now we show that $\on{Imm}_a(S^1,\mathbb R^2)$ is a manifold. 
We use the notation from the proof of \nmb!{2.4.3} with the freeness
assumption dropped.
For $c\in\on{Imm}_a(S^1,\mathbb R^2)$ we use the following mapping as the
inverse of a chart: 
\begin{align*}
C^\infty(S^1,(-\ep,\ep))\x S^1&\to \bigcup_{\th\in S^1} \mathcal Q(c(\quad+\th)) 
\East{\al}{} \on{Imm}_a(S^1,\mathbb R^2),
\\
(f,\th) &\mapsto \ps_{c(\quad+\th)}(f(\quad+\th))
  \mapsto \al(\ps_{c(\quad+\th)}(f(\quad+\th)))
\end{align*}
The chart changes are smooth: 
If for $(f_i,\th_i)\in C^\infty(S^1,(-\ep,\ep))\x S^1$ we have
$\al(\ps_{c_1(\quad+\th_1)}(f_1(\quad+\th_1)))
  =\al(\ps_{c_2(\quad+\th_2)}(f_2(\quad+\th_2)))$ 
then the initial points agree and both curves are equally oriented so that 
$c_1(\th+\th_1)+f_1(\th+\th_1)n_{c_1}(\th+\th_1)
  =c_2(\ph(\th)+\th_2)+f_2(\ph(\th)+\th_2)n_{c_2}(\ph(\th)+\th_2)$
for all $\th$. From this one can express $(f_2,\th_2)$ smoothly in terms of
$(f_1,\th_1)$.

For the latter assertion one has to show that a smooth path
through $e_1$ in $\mathcal Q(c_1)$ is mapped to a smooth path in
$\on{Diff}_1(S^1)$. This follows from the finite dimensional implicit
function theorem. 
The mapping $\al$ is now smooth into $\on{Imm}_a(S^1,\mathbb R^2)$ 
and the diffeomorphism 
$\on{Imm}(S^1,\mathbb R^2)\to\on{Imm}_{a}(S^1,\mathbb R^2)\x \on{Diff}_1(S^1)$
is given by $c\mapsto (\al(c),\si_c)$ with inverse $(e,\ph)\mapsto
e\o\ph\i$.
Only the group $S^1\ltimes \mathbb Z_2$ of rotations and reflections of $S^1$
then still acts on $\on{Imm}_a(S^1,\mathbb R^2)$ with orbit space
$B_{i}(S^1,\mathbb R^2)$.
The rest is clear.
\qed\end{demo}

\subsection*{\nmb.{2.7}. 
Tangent space, length, curvature, and Frenet-Serret formulas on 
$\on{Imm}_a(S^1,\mathbb R^2)$}
A smooth curve $t\mapsto c(\quad,t)\in \on{Imm}(S^1,\mathbb R^2)$ lies in 
$\on{Imm}_a(S^1,\mathbb R^2)$ if and only if $|\p_\th c|^2=|c_\th|^2$ is 
constant in $\th$, i.e.,
$\p_\th|c_\th|^2 = 2\langle c_\th,c_{\th\th} \rangle=0$. Thus
$h=\p_t|_0 c\in T_c\on{Imm}(S^1,\mathbb R^2)=C^\infty(S^1,\mathbb R^2)$
is tangent to $\on{Imm}_a(S^1,\mathbb R^2)$ at the foot point $c$ if and 
only if
$\langle h_\th,c_{\th\th} \rangle + 
\langle h_{\th\th},c_\th \rangle = \langle h_\th,c_\th \rangle_\th=0$,
i.e., $\langle h_\th,c_\th \rangle$ is constant in $\th$. 
For $c\in \on{Imm}_a(S^1,\mathbb R^2)$ the volume form is constant 
in $\th$ since  $|c_\th| = {\ell(c)}/{2\pi}.$ Thus for the curvature we have 
$$
\ka:\on{Imm}_a(S^1,\mathbb R^2)\to C^\infty(S^1,\mathbb R),\quad
\ka(c)= \Bigl(\frac{2\pi}{\ell(c)}\Bigr)^3\det(c_\th,c_{\th\th})
= \Bigl(\frac{2\pi}{\ell(c)}\Bigr)^3\langle ic_\th,c_{\th\th}\rangle
$$
and for the derivative of the length function we get
$$
d\ell(c)(h) = \int_{S^1} \frac{\langle h_\th,c_\th \rangle}{|c_\th|}d\th
= \frac{(2\pi)^2}{\ell(c)}\langle h_\th(1),c_\th(1) \rangle.
$$
Since $c_{\th\th}$ is orthogonal to $c_\th$ we have (Frenet formulas)
\begin{align*}
c_{\th\th} 
&= \Bigl(\frac{2\pi}{\ell(c)}\Bigr)^2\langle ic_\th,c_{\th\th} \rangle ic_\th
     = \frac{\ell(c)}{2\pi} \ka(c)ic_\th,\\
c_{\th\th\th} 
&= \frac{\ell(c)}{2\pi} \ka(c)_\th\,ic_\th
     + \frac{\ell(c)}{2\pi} \ka(c)ic_{\th\th}
= \frac{\ell(c)}{2\pi} \ka(c)_\th\,ic_\th
     - \Bigl(\frac{\ell(c)}{2\pi}\Bigr)^2 \ka(c)^2 c_{\th}.
\end{align*}
The derivative of the curvature thus becomes:
\begin{align*}
d\ka(c)(h) 
%&= -3\Bigl(\frac{2\pi}{\ell(c)}\Bigr)^5
%     \langle h_\th(1),c_\th(1) \rangle \langle ic_\th,c_{\th\th} \rangle
%     + \Bigl(\frac{2\pi}{\ell(c)}\Bigr)^3
%     \bigl(\langle ih_\th,c_{\th\th} \rangle
%     +\langle ic_\th,h_{\th\th}\rangle\bigr)\\
%&= -3\Bigl(\frac{2\pi}{\ell(c)}\Bigr)^2
%     \langle h_\th,c_\th\rangle \ka(c)
%     + \Bigl(\frac{2\pi}{\ell(c)}\Bigr)^2
%     \ka(c)\langle ih_\th,ic_{\th} \rangle
%     + \Bigl(\frac{2\pi}{\ell(c)}\Bigr)^3
%     \langle ic_\th,h_{\th\th}\rangle\bigr)\\
&= -2\Bigl(\frac{2\pi}{\ell(c)}\Bigr)^2
     \langle h_\th,c_\th\rangle \ka(c)
     + \Bigl(\frac{2\pi}{\ell(c)}\Bigr)^3
     \langle ic_\th,h_{\th\th}\rangle\bigr).\\
\end{align*}
%In particular we have 
%$$
%\ka(c)_\th = \Bigl(\frac{2\pi}{\ell(c)}\Bigr)^3
%     \langle ic_\th,c_{\th\th\th}\rangle
%$$

\subsection*{\nmb.{2.8}. Horizontality on $\on{Imm}_a(S^1,\mathbb R^2)$}
Let us denote by $\on{Imm}_{a,f}(S^1,\mathbb R^2)$ the splitting submanifold
of $\on{Imm}$ consisting of all constant speed free immersions.
 From \nmb!{2.6} and \nmb!{2.4.3} we conclude that the projection
$\on{Imm}_{a,f}(S^1,\mathbb R^2)\to B_f(S^1,\mathbb R^2)$ is principal fiber
bundle with structure group $S^1\ltimes \mathbb Z_2$, and it is a reduction of
the principal fibration $\on{Imm}_f\to B_f$. 
The principal connection described in \nmb!{2.4.3} is not compatible with
this reduction. But we can easily find some principal connections. 
The one we will use is described by the horizontal bundle with fiber
$\mathcal N_a,c$ consisting of all vector fields $h$ along $c$ such that 
$\langle h_\th,c_\th \rangle_\th=0$ (tangent to $\on{Imm}_a$) and $\langle
h(1),c_\th(1) \rangle=0$ for $1\in S^1$ (horizontality). This connection
admits a smooth parallel transport; but we can even do better, beyond the
principal bundle, in the following proposition whose 
proof is similar and simpler than that of proposition
\nmb!{2.5}.

\begin{proclaim}{Proposition}
For any smooth path $c$ in $\on{Imm}_a(S^1,\mathbb R^2)$ there exists a smooth
curve $\ph_c$ in $S^1$ with $\ph_c(0)=1$ depending smoothly on $c$ such that
the path $e$ given by $e(t,\th)=c(t,\ph_c(t)\th)$ is horizontal:
$e_t(1)\bot e_\th(1)$. \qed
\end{proclaim}

\subsection*{\nmb.{2.9}. The degree of immersions
}

Recall that the degree of an immersion $c:S^1\to \mathbb R^2$ is the winding
number with respect to 0 of the tangent $c':S^1\to \mathbb R^2$. Since this is
invariant under isotopies of immersions, the manifold 
$\on{Imm}(S^1,\mathbb R^2)$ decomposes into the disjoint union of the 
open submanifolds $\on{Imm}^k(S^1,\mathbb R^2)$ for $k\in \mathbb Z$ 
according to the degree $k$.
We shall also need the space $\on{Imm}^k_a(S^1,\mathbb R^2)$ of all immersions
of degree $k$ with constant speed.

\begin{proclaim}{\nmb.{2.10}. Theorem} 
\begin{enumerate}
\item[(\nmb:{1})]
  The manifold $\on{Imm}^k(S^1,\mathbb R^2)$ of immersed 
  curves of degree $k$ contains the subspace $\on{Imm}^k_a(S^1,\mathbb R^2)$
  as smooth strong deformation retract.
\item[(\nmb:{2})] For $k\ne 0$
  the manifold $\on{Imm}^k_a(S^1,\mathbb R^2)$ of immersed constant speed 
  curves of degree $k$ contains $S^1$ as a strong smooth deformation retract.
\item[(\nmb:{3})]
  For $k\ne 0$ the manifold $B_i^{k}(S^1,\mathbb R^2):=\on{Imm}^k(S^1,\mathbb
  R^2)/\on{Diff}^+(S^1)$ is contractible. 
\end{enumerate}
\end{proclaim}
Note that for $k\ne 0$ $\on{Imm}^k$ is invariant under the action of
the group $\on{Diff}^+(S^1)$ of orientation preserving diffeomorphism only,
and that any orientation reversing diffeomorphism maps $\on{Imm}^k$ to
$\on{Imm}^{-k}$. 

The nontrivial $S^1$ in $\on{Imm}^k$ appears in 2 ways: ($a$) by rotating 
each curve around $c(0)$ so that $c'(0)$ rotates. 
And ($b$) also by acting $S^1\ni
\be\mapsto(c(\th)\mapsto c(\be\th))$. The two corresponding elements $a$
and $b$ in the fundamental group are then related by $a^k=b$ which explains
our failure to describe the topological type of $B_i^0$.

\begin{demo}{Proof}
{\rm (\nmb!{1})} is a consequence of \nmb!{2.6} since $\on{Diff}_1^+(S^1)$ is
contractible.

The general proof is inspired by the proof of the Whitney-Graustein
theorem, \cit!{6}, \cit!{2}, \cit!{3}.
We shall view curves here as $2\pi$-periodic plane-valued functions.
For any curve $c$ we consider its {\it center of mass}
$$
C(c)=\text{Center}(c) := \frac1{\ell(c)}\int_0^{2\pi}c(u)\,|c'(u)|\,du\in \mathbb R^2
$$
which is invariant under $\on{Diff}(S^1)$. We shall also use
$\al(c)=c'(0)/|c'(0)|$.

\noindent{\bf The case $k\ne 0$.} 
We first embedd $S^1$ into $\on{Imm}(S^1,\mathbb R^2)$ in the following way.
For $\al\in S^1\subset\mathbb C=\mathbb R^2$ and $k\ne 0$ 
we but $e_\al(\th) = \al.e^{ik\th}/ik$, a circle 
of radius $1/|k|$ transversed $k$-times in the direction indicated by the 
sign of $k$. 
Note that we have $\text{Center}(e_\al)=0$ and $e_\al'(0)=\al$.

Since the isotopies to be constructed later
will destroy the property of having constant speed, 
we shall first construct a smooth deformation retraction $A:[0,1]\x \on{Imm}^k \to
\on{Imm}^k_{1,0}$ onto the subspace $\on{Imm}^k_{1,0}$ of unit speed degree
$k\ne 0$ curves with center $0$.

Let $c:\mathbb R\to \mathbb R^2$ be an arbitrary constant speed immersion of 
degree $k$, period $2\pi$, and length $\ell(c)$. 
Let $s_c(v)=\int_0^v|c'(u)|\,du$ be the arc-length function of $c$ and put
\begin{displaymath}
A(c,t,u) =\Bigl(1-t+t\frac{2\pi}{\ell(c)}\Bigr).
\Bigl(c\bigl((1-t)u+t.s_c\i(\tfrac{\ell(c)}{2\pi}u)\bigr) -t.C(c)\Bigr).
\end{displaymath}
Then $A_c$ is an isotopy between $c$ and $c_1:=A(c,1,\quad)$ depending smoothly
on $c$. The immersion  
$c_1$ has unit speed, length $2\pi$, and $\text{Center}(c_1)=0$.
Moreover, for the winding number $w_0$ around 0 we have: 
\begin{displaymath}
w_0(c_1'|_{[0,2\pi]}) = \text{deg}(c_1) = \on{deg}(c) = k =
\text{deg}(e_{\al(c)})
= w_0(e_{\al(c)}'|_{[0,2\pi]}).
\end{displaymath}
Thus $\on{Imm}^k$ contains the space $\on{Imm}^k_{1,0}$ of unit speed
immersions with center of mass 0 and degree $k$ 
as smooth strong deformation retract.

For $c\in\on{Imm}^k_{1,0}$ a unit speed immersion with center $0$ 
we now construct an isotopy $t\mapsto H^1(c,t,\quad)$ between $c$ and
a suitable curve $e_{\al}$. It will destroy the unit speed property,
however.
For $d\on{arg}=\frac{-xdy+ydx}{\sqrt{x^2+y^2}}$ we put: 
\begin{align*}
\ph_c(u)&:=\int_{c'|[0,u]}d\on{arg},
  \qquad\quad\text{ so that }
  c'(u)=c'(0)\,e^{i\ph_c(u)},
\\
\al(c)&:=\frac1{2\pi}\int_0^{2\pi}(\ph_c(v)-kv)dv,
\\
\ps_c(t,u)&:=(1-t)\ph_c(u)+t(ku+\al(c)),
\\
h(c,t,u) &:= \int_0^u e^{i\ps_c(t,v)}dv 
  - \frac{u}{2\pi}\int_0^{2\pi} e^{i\ps_c(t,v)}dv,
\\
H^1(c,t,u) &:= c'(0)\Bigl(h(c,t,u) -\text{Center}(h(c,t,\quad)\Bigr)
\end{align*}
Then $H^1(c,t,u)$ is smooth in all variables, 
$2\pi$-periodic in $u$, with center of mass at 0,
$H^1(1,c,u)$ equals one the $e_\al$'s, and $H^1(0,c,u)=c(u)$. 
But $H^1(c,t,\quad)$ is, however, no longer of unit speed in general.
And we still have to show that $t\mapsto h(c,t,\quad)$ (and consequently
$H^1$) is an isotopy.
\begin{gather}
\p_u h(c,t,u) = e^{i\ps_c(t,u)} - \frac1{2\pi}\int_0^{2\pi}
  e^{i\ps_c(t,v)}dv,
\notag\\
\Bigl|\frac1{2\pi}\int_0^{2\pi} e^{i\ps_c(t,v)}dv \Bigr| \le 1.
\AMStag{\nmb:{4}}
\end{gather}
If the last inequality is strict we have $\p_u h(t,u)\ne 0$ so that $h$ is
an isotopy. 
If we have equality then $\ps_c(t,v)$ is constant in $v$ which
leads to a contradiction as follows: If $k\ne 0$ then
$\ps_c(t,2\pi)-\ps_c(t,0)=2\pi k$ so it cannot be constant for any $t$. 

Let us finally check how this construction depends on the choice of the
base point $c(0)$. We have:
\begin{align*}
\ph_{c(\be+\quad)}(u)&=\ph_c(\be+u)-\ph_c(\be),
\\
\al(c(\be+\quad))&=\al(c)+k\be-\ph_c(\be),
\\
\ps_{c(\be+\quad)}(t,u)&=\ps_c(t,u+\be)-\ph_c(\be),
\\
h(c(\be+\quad),t,u) &= e^{-i\ph_c(\be)}(h(c,t,\be+u)-h(c,t,\be)),
\\
H^1(c(\be+\quad),t,u) &= H^1(c,t,\be+u).
\end{align*}
Let us now deform $H^1$ back into $\on{Imm}^k_{1,0}$.
For $c\in \on{Imm}^k_{1,0}$ we consider 
\begin{align*}
H^2(c,t,u) &:= A(1,H^1(c,t,\quad),u),
\\
H^3(c,t,u) &:= H^2(c,t,u+\ph_{H^2(c)}(t)),
\end{align*}
where the $\ph_f$ for a unit speed path $f$ is from proposition \nmb!{2.8}, 
so that $H^3(c)$ is a horizontal path of unit speed curves of length $2\pi$,
(i.e., $\p_tH^3(c,t,0)\bot \p_u|_0H^3(c,t,u)$). 

The isotopy $A$ reacts in a complicated way to rotations of the parameter,
but we have $A(c(\be+\quad),1,u)=A(c,1,\frac{2\pi}{\ell(c)}s_c(\be)+u)$.
Thus $H^3(c(\quad+\be),t,u)=H^3(c,t,u+\be)$, so $H^3$ is 
equivariant under the rotation group $S^1\subset \on{Diff}(S^1)$. 
For $k\ne 0$ we get 
an equivariant smooth strong deformation retract within 
$\on{Imm}^k_{1,0}$ onto the subset $\{e_\al:\al\in S^1\}\subset\on{Imm}^k_{1,0}$
which is invariant under the rotation group $S^1\subset \on{Diff}(S^1)$.
It factors to a smooth contraction on $B_i^k$. 
This proves assertions (\nmb|{2}) and (\nmb|{3})
for $k\ne 0$. 
\qed\end{demo}

\subsection*{\nmb.{2.11} Bigger spaces of `immersed' curves
}

We want to introduce a larger space containing 
$B_i(S^1, \mathbb R^2)$, which is complete in a suitable metric. This will serve
as an ambient space which will contain the completion of $B_i(S^1, \mathbb R^2)$.
Let $\on{Cont}(S^1,\mathbb R^2)$ be the space of all {\it continuous} functions
$c: S^1 \rightarrow \mathbb R^2$. Instead of a group operation and its associated
orbit space, we introduce an equivalence relation on $\on{Cont}(S^1,\mathbb R^2)$.
Define a subset $R \subset S^1 \times S^1$ to be a {\it monotone correspondence}
if it is the image of a map 
\begin{align*} 
& x \rightarrow (h(x) \mod 2\pi, k(x) \mod 2\pi), \quad \text{where} 
\\& 
h, k: \mathbb R \rightarrow \mathbb R \quad 
  \text{are monotone non-decreasing continuous functions such that} 
\\& 
h(x+2\pi) \equiv h(x)+2\pi, k(x+2\pi) \equiv k(x)+2\pi.
\end{align*}In words, this is
an orientation preserving homeomorphism from $S^1$ to $S^1$ which is allowed to
have intervals where one or the other variable remains constant while the other
continues to increase. (These correspondences arise naturally in computer
vision in comparing the images seen by the right and left eyes, see \cit!{8}.)
Then we define the equivalence relation on $\on{Cont}(S^1,\mathbb R^2)$ by
$c \sim d$ if and only if there is a monotone correspondence $R$ such that
for all $\th,\ph \in R, c(\th) = d(\ph).$ It is easily seen that any non-constant
$c \in \on{Cont}(S^1,\mathbb R^2)$ is equivalent to an $c_1$ which is not constant
on any intervals in $S^1$ and that for such $c_1$'s and $d_1$'s, the equivalence
relation amounts to $c_1 \circ h \equiv d_1$ for some homeomorphism $h$ of
$S^1$. Let $B_i^{cont}(S^1,\mathbb R^2)$ be the quotient space by this 
equivalence relation. We call these {\it Fr\'echet curves}.

The quotient metric on $B_i^{cont}(S^1,\mathbb R^2)$ is called the 
{\it Fr\'echet} metric, a variant of the {\it Hausdorff} metric mentioned in the
Introduction, both being $L^\infty$ type metrics. Namely, define
\begin{align*}
d_\infty (c,d) &= \inf_{\text{monotone corresp.} R}\quad \Bigl(
\sup_{(\th,\ph) \in R} |c(\th)-d(\ph)|\Bigr) \\
& = \inf_{\text{homeomorph.} h:S^1 \rightarrow S^1} 
\|c\circ h - d\|_\infty.
\end{align*}It is straightforward to check that this makes 
$B_i^{cont}(S^1,\mathbb R^2)$ into a complete metric space.

Another very natural space is the subset 
$B_i^{lip}(S^1,\mathbb R^2) \subset B_i^{cont}(S^1,\mathbb R^2)$
given by the non-constant {\it Lipschitz} maps $c:S^1 \rightarrow \mathbb R^2$.
The great virtue of Lipschitz maps is that their images are rectifiable curves
and thus each of them is equivalent to a map $d$ in which $\th$ is proportional
to arclength, as in the previous section. More precisely, if $c$ is Lipschitz,
then $c_\th$ exists almost everywhere and is bounded and we can reparametrize
by:
\begin{displaymath}
h(\th) = \int_0^\th |c_\th| d\th 
\biggl/ \int_0^{2\pi} |c_\th| d\th,
\end{displaymath}
obtaining an equivalent $d$ for which $|d_\th| \equiv L/2\pi$. This $d$ will
be unique up to rotations, i.e.\ the action of $S^1$ in the previous section.

This subspace of rectifiable Fr\'echet curves is the subject of a nice 
compactness theorem due to Hilbert, namely that the set of all such curves 
in a closed bounded subset of $\mathbb R^2$  and whose length is bounded 
is compact in the Fr\'echet metric. 
This can be seen as follows: we can lift all such curves to specific Lipschitz 
maps $c$ whose Lipschitz constants are bounded. This set is an equicontinuous 
set of functions by the bound on the Lipschitz constant. By the Ascoli-Arzela
theorem the topology of pointwise convergence equals then the topology of
uniform convergence on $S^1$. So this set is a closed subset in a product
of $S^1$ copies of a large ball in $\mathbb R^2$; this product is compact. The
Fr\'echet metric is coarser than the uniform metric, so our set is also
compact.

\section*{\totoc\nmb0{3}. Metrics on spaces of curves
}

\subsection*{\nmb.{3.1} Need for invariance under reparametrization
} 

The pointwise metric on the space of
immersions $\on{Imm}(S^1,\mathbb R^2)$ is given by 
\begin{displaymath}
G_c(h,k) := \int_{S^1} \langle h(\th),k(\th) \rangle d\th.
\end{displaymath}
This Riemannian metric is not invariant under reparameterizations of the
variable $\th$ and thus does not induce a sensible metric on the quotient space 
$B_{i}(S^1,\mathbb R^2)$. Indeed, it induces the zero metric since 
{\it for any two curves $C_0,C_1\in B_{i}(S^1,\mathbb R^2)$
the infimum of the arc lengths of curves in $\on{Imm}(S^1,\mathbb R^2)$ which
connect embeddings $c_0,c_1\in \on{Imm}(S^1,\mathbb R^2)$ with $\pi(c_i)=C_i$
turns out to be zero.}
To see this, take any $c_0$ in the $\on{Diff}(S^1)$-orbit over $C_0$.
Take the following variation $c(\th,t)$ of $c_0$: for $\th$ outside a small 
neighborhood $U$ of length $\ep$ of 1 in $S^1$, $c(\th,t)=c_0(\th)$. If
$\th\in U$, then the variation for $t\in[0,1/2]$ 
moves the small part of $c_0$ so that $c(\th,1/2)$ for $\th$ in $U$
takes off $C_0$, goes to $C_1$, traverses nearly all of $C_1$, and returns
to $C_0$. Now in the orbit through $c(\cdot,1/2)$, reparameterize in such a
way that the new curve is diligently traversing $C_1$ for $\th\notin U$,
and for $\th\in U$ it travels back to $C_0$, runs along $C_0$, and comes
back to $C_1$. This reparametrized curve is then varied for $t\in
[1/2,1]$ in such a way, that the part for $\th\in U$ is moved towards
$C_2$. It is clear that the length of both variations is bounded by a constant
(depending on the distance between $C_0$ and $C_1$ and the lengths of both
$C_0$ and $C_1$) times $\ep$.

\subsection*{\nmb.{3.2}. The simplest Riemannian metric on $B_i$}

Let $h,k\in C^\infty(S^1,\mathbb R^2)$ be two tangent vectors with foot 
point $c\in \on{Imm}(S^1,\mathbb R^2)$. 
The induced volume form is
$\on{vol}(c)=\langle \p_\th c,\p_\th c\rangle^{1/2}\,d\th = |c_\th|d\th$. 
We consider first the simple $H^0$ weak Riemannian metric on 
$\on{Imm}(S^1,\mathbb R^2)$:
\begin{equation}
G_c(h,k) := \int_{S^1}\langle h(\th),k(\th) \rangle |c'(\th)|\,d\th
\AMStag{\nmb:{1}}\end{equation}
which is invariant under $\on{Diff}(S^1)$. This makes the map 
$\pi:\on{Imm}(S^1,\mathbb R^2) \rightarrow B_i(S^1,\mathbb R^2)$ into a
{\it Riemannian submersion} (off the singularities of $B_i(S^1,\mathbb R^2)$) 
which is very convenient. We call this
the \idx{\it $H^0$-metric}.

Now we can determine the bundle $\mathcal N\to \on{Imm}(S^1,\mathbb R^2)$ 
of tangent vectors which are normal to the $\on{Diff}(S^1)$-orbits.
The tangent vectors to the orbits are $T_c(c\o\on{Diff}(S^1)) = 
\{g.c_\th:g\in C^\infty(S^1,\mathbb R)\}$. Inserting this for $k$
into the expression \thetag{\nmb|{1}} of the metric we see that 
\begin{align}
\mathcal N_c & =\{h\in C^\infty(S^1,\mathbb R^2): \langle h,c_\th\rangle =0\}
     \AMStag{\nmb:{2}}\\
&=\{aic_\th\in C^\infty(S^1,\mathbb R^2): a\in C^\infty(S^1,\mathbb R)\}
\notag\\
&=\{bn_c\in C^\infty(S^1,\mathbb R^2): b\in C^\infty(S^1,\mathbb R)\},
\notag\end{align}
where $n_c$ is the normal unit field along $c$. 

A tangent vector $h\in T_c\on{Imm}(S^1,\mathbb R^2)=C^\infty(S^1,\mathbb R^2)$ 
has an orthonormal decomposition 
\begin{align}
h&=h^\top+h^\bot \in T_c(c\o\on{Diff}^+(S^1)) \oplus \mathcal N_c\qquad
  \text{ where}
\notag\\
h^\top &= \frac{\langle h,c_\th \rangle}{|c_\th|^2}c_\th 
  \in  T_c(c\o\on{Diff}^+(S^1)),
\AMStag{\nmb:{3}}
\\
h^\bot &= \frac{\langle h,ic_\th \rangle}{|c_\th|^2}ic_\th \in  \mathcal N_c,
\notag\end{align}
into smooth tangential and normal components.

Since the Riemannian metric $G$ on $\on{Imm}(S^1, \mathbb R^2)$ is invariant
under the action of $\on{Diff}(S^1)$ it induces a metric on the quotient
$B_i(S^1, \mathbb R^2)$ as follows. For any $C_0, C_1 \in B_i$, consider all
liftings $c_0, c_1 \in \on{Imm}$ such that $\pi(c_0)=C_0, \pi(c_1)=C_1$ and
all smooth curves $t\mapsto (\th\mapsto c(t,\th))$ in $\on{Imm}(S^1,\mathbb R^2)$
with $c(0,\cdot)=c_0$ and $c(1,\cdot)=c_1$. 
Since the metric $G$ is invariant under the action of $\on{Diff}(S^1)$
the arc-length of the curve 
% $t\mapsto \pi(c(t,\cdot))$ in $B_{i,f}(S^1,\mathbb R^2)$ is given by 
$t\mapsto \pi(c(t,\cdot))$ in $B_i(S^1,\mathbb R^2)$ is given by 
\begin{align}
L^{\text{hor}}_G(c) &:= L_G(\pi(c(t,\cdot)))
= \int_0^1 \sqrt{G_{\pi(c)}(T_c\pi.c_t,T_c\pi.c_t)}\,dt
  = \int_0^1 \sqrt{G_{c}(c_t^\bot,c_t^\bot)}\,dt
\notag\\&
=\int_0^1\Bigl(\int_{S^1}
\left\langle \frac{\langle c_t,ic_\th \rangle}{|c_\th|^2}ic_\th, 
\frac{\langle c_t,ic_\th \rangle}{|c_\th|^2}ic_\th \right\rangle
|c_\th|\,d\th\Bigr)^{\tfrac12}dt
\notag\\&
=\int_0^1\Bigl(\int_{S^1}
\langle c_t,n_c \rangle^2 |c_\th|\,d\th \Bigr)^{\tfrac12}dt
\AMStag{\nmb:{4}}
\\&
=\int_0^1\Bigl(\int_{S^1}
\langle c_t,ic_\th \rangle^2 \frac{d\th}{|c_\th|} \Bigr)^{\tfrac12}dt
\notag\end{align}
The metric on $B_i(S^1, \mathbb R^2)$ is defined by taking the infimum of this 
over all paths $c$ (and all lifts $c_0, c_1)$:
\begin{displaymath}
\on{dist}^{B_i}_{G}(C_1,C_2) = 
\inf_c L^{\text{hor}}_G(c).
\end{displaymath}
%
%Note that if a path $\pi(c)$ in $B_i(S^1, \mathbb R^2)$ is given, then one can 
%choose its lift to a path $c$ in $\on{Imm}(S^1, \mathbb R^2)$ to have various
%good properties. Firstly, we can choose the lift $c(0,\quad)$ 
%of the inital curve to
%have a parametrization of constant speed, i.e.\ if its length is $\ell$, then
%$|c_\th|(\th,0) = \ell/2\pi$ for all $\th \in S^1$. Secondly, we can make the 
%tangent vector to $c$ everywhere horizontal, i.e.\ $<c_t, c_\th> \equiv 0$,
%by \nmb!{2.5}.
%Thirdly, we can reparametrize the coordinate $t$ on the path of length $L$
%so that the path is traversed at constant speed, i.e.\ 
%$$\int_{S^1} <c_t,ic_\th>^2 d\th/|c_\th| \equiv L^2, 
%\text{ for all }
%0 \le t \le 1.$$
%
Unfortunately, we will see below that this metric is too weak: the distance 
that it defines turns out to be identically zero! For this reason, we will 
mostly study in this paper a family of stronger metrics. These are obtained 
by the most minimal change in $G$. We want to preserve two simple properties 
of the metric: that it is local and that it has no derivatives in it. The 
standard way to strengthen the metric is go from an $H^0$ metric to an $H^1$ 
metric. But when we work out the natural $H^1$ metric, picking out those terms 
which are local and do not involve derivatives leads us to our chosen metric.

We consider next the $H^1$ weak Riemannian metric on 
$\on{Imm}(S^1,\mathbb R^2)$:
\begin{equation}
G^1_c(h,k) := \int_{S^1}
\bigl(\langle h(\th),k(\th)\rangle + A \frac{\langle h_\th,k_\th \rangle}
{|c_\th|^2}
\bigr) |c_\th|\,d\th.
\AMStag{\nmb:{5}}\end{equation}
which is invariant under $\on{Diff}(S^1)$. Thus 
$\pi:\on{Imm}(S^1,\mathbb R^2) \rightarrow B_i(S^1,\mathbb
R^2)$ is again a
{\it Riemannian submersion} off the singularities of
$B_i(S^1,\mathbb R^2)$.
We call this the \idx{\it $H^1$-metric} on $B_i$.

To understand this metric better, we assume $h=k=a\frac{ic_\th}{|c_\th|}
+b\frac{c_\th}{|c_\th|}$. Moreover, for any function $f(\th)$, we write
$f_s = \frac{f_\th}{|c_\th|}$ for the derivative with respect
to arc length. Then:
$$h_s= \frac{h_\th}{|c_\th|}=(aic_s+bc_s)_s=(a_s+\ka b)ic_s+(b_s-\ka a)c_s.$$
Therefore:
\begin{align*} 
G^1_c(h,h) &= \int_{S^1} \bigl( a^2+b^2+A(a_s+\ka b)^2+A(b_s-\ka a)^2 
\bigr) ds\\
&= \int_{S^1} \bigl(a^2(1+A\ka^2)+Aa^2_s \bigr) + 2A\ka (a_s b - b_s a) +
\bigl( b^2(1+A\ka^2) + Ab^2_s \bigr) ds
\end{align*}
Letting $T_1$ and $T_2$ be the differential operators 
$T_1 = I+A\ka^2 - A(\frac{d}{ds})^2$, $T_2 = A(\ka_s + 2\ka \frac{d}{ds})$, 
then integrating by parts on $S^1$, we get:
$$ G^1_c(h,h) = \int_{S^1} \bigl( T_1(a).a + 2 T_2(a).b + T_1(b).b \bigr) ds.$$
Note that $T_1$ is a positive definite self-adjoint operator on functions on $c$, 
hence it has an inverse given by a Green's function which we write $T_1^{-1}$.
Completing the square and using that $T_1$ is self-adjoint, 
we simplify the metric to:
$$ G^1_c(h,h)\! =\! \int_c \Bigr( T_1(a).a - T_1^{-1}(T_2(a)).T_2(a) + 
T_1\bigl(b+T_1^{-1}(T_2(a))\bigr).\bigl(b+T_1^{-1}(T_2(a))\bigr) \Bigl) ds.$$

If we fix $a$ and minimize this in $b$, we get the bundle 
$\mathcal N^1\to \on{Imm}(S^1,\mathbb R^2)$ 
of tangent vectors which are $G^1$-normal to the $\on{Diff}(S^1)$-orbits. In
other words: 
\begin{align*}
\mathcal N^1_c & =\{h\in C^\infty(S^1,\mathbb R^2): h = aic_s+bc_s, 
b=-T_1^{-1}(T_2(a)) \}
\end{align*}
and on horizontal vectors of this type:
$$ G^1_c(h,h) = \int_c \bigr((1+A\ka^2)a^2 + Aa^2_s \bigl)ds -
\int_c T_1^{-1}(T_2(a)).T_2(a) ds.$$
If we drop terms involving $a_s$, say because we assume $|a_s|$ is small,
then what remains is just the integral of $(1+A\ka^2)a^2$ plus the integral
of $T_1^{-1}(\ka_s a)\ka_s a$. The second is a non-local regular integral 
operator, so dropping this we are left with the main metric of this paper:
$$ G^A_c(h,h) = \int_c (1+A\ka^2) a^2 ds, h = a ic_s$$
which we call the \idx{\it $H^0_\ka$-metric with curvature weight $A$}.
For further reference, 
on $\on{Imm}(S^1,\mathbb R^2)$, for a constant $A\ge 0$, it is
given by
\begin{equation}
G^A_c(h,k) := \int_{S^1}(1+A\ka_c(\th)^2)\langle h(\th),k(\th) \rangle 
     |c'(\th)|\,d\th
\AMStag{\nmb:{6}}\end{equation}
which is again invariant under $\on{Diff}(S^1)$. 
Thus 
$\pi:\on{Imm}(S^1,\mathbb R^2) \rightarrow B_i(S^1,\mathbb
R^2)$ is again a {\it Riemannian submersion} off the
singularities.
Note that for this metric \thetag{\nmb|{6}}, the bundle 
$\mathcal N\subset T\on{Imm}(S^1,\mathbb R^2)$ is the same as for $A=0$,
as described in \thetag{\nmb|{2}}.
The arc-length of a curve 
$t\mapsto \pi(c(t,\cdot))$ in $B_i(S^1,\mathbb R^2)$ is given
by the analogon of \thetag{\nmb|{4}}
\begin{align}
L^{\text{hor}}_{G^A}(c) &:= L_{G^A}(\pi(c(t,\cdot)))
= \int_0^1 \sqrt{G^A_{\pi(c)}(T_c\pi.c_t,T_c\pi.c_t)}\,dt
  = \int_0^1 \sqrt{G^A_{c}(c_t^\bot,c_t^\bot)}\,dt
\notag\\&
=\int_0^1\Bigl(\int_{S^1}(1+A\ka_{c}^2)
\langle c_t,n_c \rangle^2 |c_\th|\,d\th \Bigr)^{\tfrac12}dt
\AMStag{\nmb:{7}}
\\&
=\int_0^1\Bigl(\int_{S^1}(1+A\ka_{c}^2)
\langle c_t,ic_\th \rangle^2 \frac{d\th}{|c_\th|} \Bigr)^{\tfrac12}dt
\notag\end{align}
The metric on $B_i(S^1, \mathbb R^2)$ is defined by taking the infimum of this 
over all paths $c$ (and all lifts $c_0, c_1)$:
\begin{displaymath}
\on{dist}^{B_i}_{G^A}(C_1,C_2) = 
\inf_c L^{\text{hor}}_{G^A}(c).
\end{displaymath}
Note that if a path $\pi(c)$ in $B_i(S^1, \mathbb R^2)$ is given, then one can 
choose its lift to a path $c$ in $\on{Imm}(S^1, \mathbb R^2)$ to have various
good properties. Firstly, we can choose the lift $c(0,\quad)$ 
of the inital curve to
have a parametrization of constant speed, i.e.\ if its length is $\ell$, then
$|c_\th|(\th,0) = \ell/2\pi$ for all $\th \in S^1$. Secondly, we can make the 
tangent vector to $c$ everywhere horizontal, i.e.\ $<c_t, c_\th> \equiv 0$,
by \nmb!{2.5}.
Thirdly, we can reparametrize the coordinate $t$ on the path of length $L$
so that the path is traversed at constant speed, i.e.\ 
$$\int_{S^1} (1+A\kappa_c^2) <c_t,ic_\th>^2 d\th/|c_\th| \equiv L^2, 
\text{ for all }
0 \le t \le 1.$$

\subsection*{\nmb.{3.3}. A Lipschitz bound for arc length in $G^A$}

We apply the Cauchy-Schwarz inequality to the 
derivative \nmb!{2.2.4} of the length function along a path
$t\mapsto c(t,\quad)$:
\begin{align*}
\p_t\ell(c)&=d\ell(c)(c_t)
=-\int_{S^1} \ka(c)\langle c_t, n_c\rangle|c_\th| d\th
%\\&
\le     \Bigl|\int_{S^1} \ka(c)\langle c_t, n_c\rangle|c_\th| d\th\Bigr|
\\&
\le \Bigl(\int_{S^1}1^2 |c_\th|d\th\Bigr)^{\frac12}
 \Bigl(\int_{S^1} \ka(c)^2\langle c_t, n_c\rangle^2|c_\th| d\th \Bigr)
 ^{\frac12}
\\&
\le \ell(c)^{\frac12}\frac1{\sqrt{A}}
 \Bigl(\int_{S^1}(1+A \ka(c)^2)\langle c_t, n_c\rangle^2|c_\th| d\th \Bigr)
 ^{\frac12}
\end{align*}
Thus 
\begin{displaymath}
\p_t(\sqrt{\ell(c)})=\frac{\p_t\ell(c)}{2\sqrt{\ell(c)}}
\le \frac1{2\sqrt{A}}
 \Bigl(\int_{S^1}(1+A \ka(c)^2)\langle c_t, n_c\rangle^2|c_\th| d\th \Bigr)^{\frac12}
\end{displaymath}
and by using \thetag{\nmb!{3.2.7}} we get
\begin{align}
\sqrt{\ell(c_1)}-\sqrt{\ell(c_0)}
&=\int_0^1\p_t(\sqrt{\ell(c)})\,dt
\notag\\&
\le \frac1{2\sqrt{A}}
 \int_0^1\Bigl(\int_{S^1}(1+A \ka(c)^2)
        \langle c_t, n_c\rangle^2|c_\th| d\th \Bigr)^{\frac12}dt
\notag\\&
= \frac1{2\sqrt{A}}L^{\text{hor}}_{G^A}(c).
\AMStag{\nmb:{1}}
\end{align}
If we take the infimum over all paths connecting $c_0$ with the
$\on{Diff}(S^1)$-orbit through $c_1$ we get:

\begin{proclaim}{Lipschitz continuity of 
$\sqrt{\ell}:B_i(S^1,\mathbb R^2)\to \mathbb R_{\ge0}$}
For $C_0$ and $C_1$ in $B_i(S^1,\mathbb R^2)=\on{Imm}(S^1,\mathbb
R^2)/\on{Diff}(S^1)$ we have for $A>0$: 
\begin{equation}
\sqrt{\ell(C_1)}-\sqrt{\ell(C_0)}\le
\frac1{2\sqrt{A}}\on{dist}^{B_i}_{G^A}(C_1,C_2).
\AMStag{\nmb:{2}}
\end{equation}
\end{proclaim}

\subsection*{\nmb.{3.4}. Bounding the area swept by a path in $B_i$}

Secondly, we want to bound the area swept out by a path 
starting from $C_0$ to reach any
curve $C_1$ nearby in our metric. First we use the Cauchy-Schwarz inequality 
in the Hilbert space $L^2(S^1,|c_\th(t,\th)|d\th)$ to get  
\begin{multline*}
\int_{S^1}1.|c_t(t,\th)\|c_\th(t,\th)|d\th = \langle 1,|c_t|\rangle_{L^2}\le
\\
\le \|1\|_{L^2}\|c_t\|_{L^2}
= \Bigl(\int_{S^1}|c_\th(t,\th)|\,d\th\Bigr)^{\tfrac12}
\Bigl(\int_{S^1}|c_t(t,\th)|^2|c_\th(t,\th)|\,d\th\Bigr)^{\tfrac12}.
\end{multline*}
Now we assume that the variation $c(t,\th)$ is horizontal, so that $\langle
c_t,c_\th \rangle=0$. Then $L_{G^A}(c)=L^{\text{hor}}_{G^A}(c)$. 
We use this inequality and then the intermediate value theorem of integral
calculus to obtain 
\begin{align*}
L^{\text{hor}}_{G^A}(c) &= L_{G^A}(c) = \int_0^1 \sqrt{G^A_{c}(c_t,c_t)}\,dt
\\&
=\int_0^1\Bigl(\int_{S^1}(1+A\ka(c)^2)|c_t(t,\th)|^2|c_\th(t,\th)|\,d\th\Bigr)^{\tfrac12}dt
\\&
\ge\int_0^1\Bigl(\int_{S^1}|c_t(t,\th)|^2|c_\th(t,\th)|\,d\th\Bigr)^{\tfrac12}dt
\\&
\ge\int_0^1 \Bigl(\int_{S^1}|c_\th(t,\th)|\,d\th\Bigr)^{-\tfrac12}
\int_{S^1}|c_t(t,\th)\|c_\th(t,\th)|\,d\th\,dt
%\allowdisplaybreak
\\&
= \Bigl(\int_{S^1}|c_\th(t_0,\th)|\,d\th\Bigr)^{-\tfrac12}
\int_0^1\int_{S^1}|c_t(t,\th)\|c_\th(t,\th)|\,d\th\,dt
\\&\qquad\qquad\text{  for some intermediate value $0\le t_0\le 1$, }
\\&
= \frac{1}{\sqrt{\ell(c(t_0,\cdot))}}
\int_{[0,1]\x S^1}|\det dc(t,\th)|\,d\th\,dt.
\end{align*}

\begin{proclaim}{Area swept out bound} If $c$ is any path from $C_0$ to $C_1$, then
\begin{equation}
\begin{pmatrix} \text{area of the region swept} \\ \text{out by the variation $c$} 
  \end{pmatrix} \le \max_{t} \sqrt{\ell(c(t,\cdot))} \cdot
L^{\text{hor}}_{G^A}(c). 
\AMStag{\nmb:{1}}
\end{equation}
\end{proclaim}

This result enables us to compare the double cover
$B_i^{\text{or}}(S^1,\mathbb R^2)$ of our metric space $B_i(S^,\mathbb R^2)$
consisting of oriented unparametrized curves to the fundamental
space of geometric measure theory. Note that there is a map $h_1$ from
$B_i^{\text{or}}$ to
the space of 1-currents $\mathcal D'_1$ given by:
\begin{displaymath}
<h_1(c \text{ mod Diff}^+(S^1)),\omega> = \int_{S^1} c^*\om, 
\quad c \in \on{Imm}(S^1, \mathbb R^2).\end{displaymath}
The image $h_1(C)$ is, in fact, closed. For any $C$, define the integer-valued 
measurable function $w_C$ on $\mathbb R^2$ by:
\begin{displaymath}
w_C((x,y)) = \text{winding number of } C \text{ around } (x,y).\end{displaymath}
Then it is easy to see that, as currents, $h_1(C) = \partial (w_C dxdy)$, 
hence $\partial h_1(C) = 0$.

Although $h_1$ is obviously injective on the space $B_e$, 
it is not injective on $B_i$ as illustrated in Figure 1 below. The 
image of this mapping lies in the basic subset 
$\mathcal I_{1,c} \subset \mathcal D'_1$ of closed {\it integral} currents, namely those
which are both closed and countable sums of currents defined by Lipschitz
mappings $c_i: [0,1] \rightarrow \mathbb R^2$ of finite total length.  
Integral currents carry what is called the {\it flat} metric, which, for closed 
1-currents, reduces (by the isoperimetric inequality) to the area
distance 
\begin{equation}
d^\flat(C_1,C_2) = \iint_{\mathbb R^2} |w_{C_1} - w_{C_2}| dx dy.
\AMStag{\nmb:{2}}
\end{equation}

\begin{figure} \begin{center}
\epsfig{width=3in,file=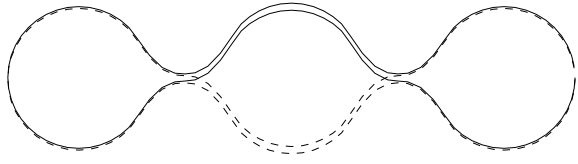}
\end{center} 
\caption{Two distinct immersions of $S^1$ in the plane whose underlying 
currents are equal. One curve is solid, the other dashed.}
\end{figure}

To connect this with our `area swept out bound', note that if we have any
path $c$ in $\on{Imm}(S^1, \mathbb R^2)$ joining $C_1$ and $C_2$, this path
defines a 2-current $w(c)$ such that $\partial w(c) = h_1(C_1)-h_1(C_2)$
and 
\begin{displaymath}
\int_{\mathbb R^2} |w(c)|\, dxdy \le \int_0^1 \int_{S^1} |\det c| d\th dt
\end{displaymath}
which is what we are calling the area swept out. But 
$\partial(w_{C_1}-w_{C_2}) = h_1(C_1)-h_1(C_2)$ too, so 
$w(c) = w_{C_1}-w_{C_2}$.
Thus
\begin{equation}
d^\flat(C_1,C_2) \le \min_{\text{all paths $c$ joining }C_1, C_2} 
\bigl[\text{area swept out by }c\bigr]
\AMStag{\nmb:{3}}
\end{equation}

Finally, we recall the fundamental compactness result of geometric measure
theory in this simple case: the space of integral 1-currents of bounded
length is compact in the flat metric. This implies that our `area swept out
bound' above has the Corollary:

\begin{proclaim}{Corollary} 
\begin{enumerate}
\item[(\nmb:{4})]
If $\{C_n\}$ is any Cauchy sequence in $B_i$ for the
metric $\on{dist}_{G^A}$, then $\{h_1(C_n)\}$ is a Cauchy sequence in 
$\mathcal I_{1,c}$ on which length is bounded.
\item[(\nmb:{5})]
Hence $h_1$ extends to a continuous map from the completion 
$\overline{B_i}$ of $B_i$ in the metric $G^A$ to $\mathcal I_{1,c}$.
\end{enumerate}
\end{proclaim}

\subsection*{\nmb.{3.5}. Bounding how far curves move in small paths in $B_i$
}

We want to bound the maximum distance a curve $C_0$ can move on
any path whose length is small in $G^A$ metric. Fix the initial curve $C_0$
and let $\ell$ be its length. The result is:

\begin{proclaim}{Maximum distance bound}
Let $\epsilon < \min\{ 2\sqrt{A\ell}, \ell^{3/2}\}/8$ and consider   
$\eta = 4(\ell^{3/4}A^{-1/4} + \ell^{1/4})\sqrt{\epsilon}$. 
Then for any path $c$ starting at $C_0$ whose length is $\epsilon$, 
the final curve lies in the tubular
neighborhood of $C_0$ of width $\eta$. More precisely, if we choose the path
$c(t,\th)$ to be horizontal, then $\max_\th |c(0,\th)-c(1,\th)| < \eta$.
\end{proclaim}

\begin{demo}{Proof}
For all of this proof, we assume the path in $B_i$ has been lifted to 
a horizontal path $c \in \on{Imm}(S^1,\mathbb R^2)$ with  
$|c_\th|(\th,0) \equiv \ell/2\pi$, so that $\langle c_t,c_\th\rangle\equiv 0$, 
and also $\int_{S^1} (1+A\kappa_c^2) |c_t|^2 |c_\th|\, d\th \equiv \epsilon^2$.
The first step in the proof is to refine the Lipschitz bound on the length of
a curve to a local estimate. Note that by horizontality
\begin{align*}
\frac{\partial}{\partial t} \sqrt{|c_\th|} 
&= \frac{<c_{\th t}, c_\th>}{2|c_\th|^{3/2}} 
%\\&
= - \frac{<c_t, c_{\th \th}>}{2|c_\th|^{3/2}} 
%\\&
= - \frac{<c_t,ic_\th>}{2|c_\th|} \kappa_c\, |c_\th|^{1/2}
%\\&
= \mp \tfrac12  \kappa_c\,|c_t| |c_\th|^{1/2}
\end{align*} hence
\begin{displaymath}
\int_{S^1} \bigl(\frac{\partial}{\partial t} \sqrt{|c_\th|} \bigr)^2 ds
\le \frac{\epsilon^2}{4A}.\end{displaymath}
Now we make the key definition:
\begin{displaymath}
\widetilde{|c_\th|}(t,\th) = \min_{0 \le t_1 \le t} |c_\th|(t_1,\th).
\end{displaymath}
Note that the $t$-derivative of $\widetilde{|c_\th|}$ is either 0 or equal to
that of $|c_\th|$ and is $\le 0$. Thus we have:
\begin{align*}
\int_{S^1} \Bigl( \sqrt{\frac{\ell}{2\pi}} 
  - \sqrt{\widetilde{|c_\th|}(1,\th)} \Bigr) d\th
&\le \int_0^1 \int_{S^1} -\frac{\partial}{\partial t} 
\sqrt{\widetilde{|c_\th|}} d\th dt \\ 
&\le \int_0^1 \int_{S^1} \Bigl| \frac{\partial}{\partial t} 
\sqrt{|c_\th|} \Bigr| d\th dt \\
&\le \int_0^1 \Bigl(\int_{S^1} d\th \Bigr)^{1/2} \cdot \Bigl( \int_{S^1} 
\Bigl| \frac{\partial}{\partial t} \sqrt{|c_\th|} \Bigr|^2 d\th \Bigr)^{1/2} dt\\
&\le \sqrt{2\pi} \cdot \frac{\epsilon}{2\sqrt{A}}.
\end{align*}
To make use of this inequality, let 
$E = \bigl\{ \th : \widetilde{|c_\th|}(1,\th) \le
\bigl(1- (A\ell)^{-1/4}\sqrt{\ep} \bigr)\ell/2\pi \bigr\}.$ 
Our assumption on $\ep$ gives $(A\ell)^{-1/4} \sqrt{\ep} < 1/2$, hence on
$S^1\setminus E$ we have $\widetilde{|c_\th|} > \ell/4\pi$.
On $E$ we have also 
$(\widetilde{|c_\th|})^{1/2} \le (1-(A\ell)^{-1/4} \sqrt{\ep}/2) \sqrt{\ell/2\pi}$.
Combining this with the previous inequality, we get (where $\mu(E)$is the
measure of $E)$:
\begin{displaymath}
\mu(E) \frac{1}{2\sqrt{2\pi}} (\frac{\ell}{A})^{1/4} \sqrt{\ep} \le 
\sqrt{2\pi}  \cdot \frac{\ep}{2\sqrt{A}}, \quad \text{hence}\quad
\mu(E) \le 2\pi \frac{\sqrt{\ep}}{(A\ell)^{1/4}}<\pi.
\end{displaymath}

We now use the lower bound on $|c_\th|$ on $S^1 - E$ 
to control $c(1,\th)-c(0,\th)$:
\begin{align*}
\int_{S^1 - E} |c(1,\th)-c(0,\th)|\, d\th &\le 
\int_0^1 \int_{S^1 - E} |c_t|\, d\th dt 
\\&
\le \sqrt{2\pi} \cdot \int_0^1 
  \Bigl( \int_{S^1-E} |c_t|^2 d\th \Bigr)^{1/2} dt 
\\& 
\le \frac{\sqrt{2\pi}}{\sqrt{\frac{\ell}{4\pi}}} 
  \int_0^1 \Bigl( \int_{S^1-E} |c_t|^2 |c_\th| d\th \Bigr)^{1/2} dt 
% \\&
\le \frac{2 \sqrt{2} \pi}{\sqrt{\ell}} \cdot \ep 
\end{align*} 
Again, introduce a small exceptional set 
$F = \bigl\{\th\; \bigl|\; \th \notin E \text{ and }
  |c(1,\th)-c(0,\th)| \ge \ell^{1/4} \sqrt{\ep} \bigr\}$.
By the inequality above, we get:
\begin{align*}
& \mu(F) \cdot \ell^{1/4} \sqrt{\ep} \le \frac{2 \sqrt{2} \pi \ep}{\sqrt{\ell}}, 
\quad \text{hence } 
 \mu(F) \le \frac{2 \sqrt{2} \pi \sqrt{\ep}}{\ell^{3/4}}<\pi.
\end{align*} The last inequality follows from the second assumption on $\ep$.
Knowing $\mu(E)$ and $\mu(F)$ gives us the lengths $|c(0,E)|$ and $|c(0,F)|$
in $\mathbb R^2$. But we need the lengths $|c(1,E)|$ and $|c(1,F)|$ too. We get
these using the fact that the whole length of $C_1$ can't be too large, by
\nmb!{3.3}:
\begin{align*}
\sqrt{|C_1|} &\le \sqrt{\ell} + \frac{\ep}{2\sqrt{A}}, \quad \text{hence} 
\\
|C_1| & \le \ell + 2\ep \sqrt{\frac{\ell}{A}} \le \ell 
  + \sqrt{\ep} \cdot \frac{\ell^{3/4}}{A^{1/4}}.
\end{align*} 
On $S^1\setminus E$ we have
$\widetilde{|c_\th|}>(1-(A\ell)^{-1/4}\sqrt{\ep})\ell/2\pi$, thus we get
\begin{align*}
|c(1,E\cup F)| &= |C_1| - |c(1,S^1\setminus (E\cup F))| \\
& \le \ell + \sqrt{\ep} \cdot \frac{\ell^{3/4}}{A^{1/4}} -
\Bigl( 1 - \frac{\sqrt{\ep}}{(A\ell)^{1/4}} \Bigr)\frac{\ell}{2\pi}
\Bigl(2\pi - \mu(E\cup F) \Bigr) \\
& \le \sqrt{\ep} \cdot \Bigl( 3\frac{\ell^{3/4}}{A^{1/4}} 
  + \sqrt{2} \ell^{1/4} \Bigr)
\end{align*} 
Finally, we can get from $c(0,\th)$ to $c(1,\th)$ by going via $c(0,\th')$ and
$c(1,\th')$ where $\th' \in S^1\setminus (E \cup F)\ne \emptyset$. Thus
\begin{align*}
 \max_\th |c(0,\th)-c(1,\th)| &\le |c(0,E\cup F)| + \ell^{1/4} \sqrt{\ep} + 
|c(E\cup F,1)| \\
& \le 4 (\ell^{3/4} A^{-1/4} + \ell^{1/4})\sqrt{\ep} \qed
\end{align*} \end{demo}

Combining this bound with the Lipschitz continuity of the square root 
of arc length, we get: 

\begin{proclaim}{\nmb.{3.6}. Corollary}
For any $A>0$, the map from $B_i(S^1,\mathbb R^2)$ in the $\on{dist}_{G^A}$ 
metric to the space $B^{cont}_i(S^1,\mathbb R^2)$ in the Fr\'echet metric 
is continuous, and, in fact,
uniformly continuous on every subset where the length $\ell$ is bounded. In 
particular, $\on{dist}_{G^A}$ is a separating metric on $B_i(S^1,\mathbb R^2)$. 
Moreover, the completion $\overline{B_i}(S^1,\mathbb R^2)$ of 
$B_i(S^1,\mathbb R^2)$ in this metric
can be identified with a subset of $B^{lip}_i(S^1,\mathbb R^2)$.
\end{proclaim}

If we iterate this bound, then we get the following:

\begin{proclaim}{\nmb.{3.7}. Corollary}
Consider all paths in $B_i$ joining curves $C_0$ and $C_1$. Let $L$ be
the length of such a path in the $\on{dist}_{G^A}$ metric and let 
$\ell_{\min}, \ell_{\max}$ be the minimum and maximum of the arc lengths 
of the curves in this path. 
Then there are parametrizations $c_0, c_1$ of $C_0$ and $C_1$ such that:
\begin{align*}
\max_\th |c_0(\th) - c_1(\th)| &\le 50 \max(L F^*, \sqrt{\ell_{\max} L F^*}), 
\mbox{ where}\\
F^* &= \max\Bigl(\frac{1}{\sqrt{\ell_{\min}}}, \sqrt{\frac{\ell_{\max}}{A}}\Bigr).
\end{align*}
\end{proclaim}

To prove this, you need only break up the path into a minimum number of
pieces for which the maximum distance bound \nmb!{3.5} holds and add together the
estimates for each piece. We will only sketch this proof which is 
straightforward. The constant 50 is just what comes out without attempting
to optimize the bound. The second option for bound, $50 \sqrt{\ell_{\max} L F^*}$ 
is just a rephrasing of the bound already in the theorem for short paths. If
the path is too long to satisfy the condition of the theorem, we break the
path at intermediate curves $C_i$ of length $\ell_i$ such that each begins
a subpath with length $\ep_i = \min(\sqrt{A\ell_i}, \ell_i^{3/2})/8$ 
and which don't overlap for more than 2:1. Thus $\sum_i \ep_i \le 2L$. 
Then apply the maximum distance bound \nmb!{3.5} to each piece, 
letting $\eta_i$ be the bound on how
far points move in this subpath {\it or any parts thereof} and verify:
$$ \eta_i \le 2\sqrt{2} \ell_i \le 16 \sqrt{2} \ep_i F^* ,$$
from which we get what we need by summing over $i$.

\subsection*{\nmb.{3.8}}
A final Corollary shows that if we parametrize any path appropriately, we get
explicit equicontinuous continuity bounds on the parametrization depending
only on $L,\ell_{\max}$ and $\ell_{\min}$. This is a step towards establishing
the existence of weak geodesics. The idea is this: instead of the horizontal 
parametrization $\langle c_t, c_\th \rangle \equiv 0$, we parametrize each
curve at constant speed $|c_\th| \equiv \ell(t)/2\pi$ where $\ell(t)$ is the
length of the $t^{th}$ curve and ask only that $\langle c_t, 
c_\th \rangle(0,t) \equiv 0$ for some base point $0\in [0,2\pi]$, see
\nmb!{2.8}. Then we get:

\begin{proclaim}{Corollary}
If a path $c(t,\th), 0 \le t \le 1$ satisfies 
\begin{align*}
&|c_\th(\th,t)| \equiv \ell(t)/2\pi  \qquad  \mbox{ for all } \th,t \\
&\langle c_t, c_\th \rangle(0,t) \equiv 0 \quad \qquad \mbox{ for all } t \mbox{ and}\\
&\int_{C_t} (1+A\ka_{C_t}^2) |\langle c_t, ic_\th \rangle |^2 d\th/|c_\th| 
\equiv L^2 \mbox{ for all } t,
\end{align*} 
then 
$$ |c(t_1,\th_1)-c(t_2,\th_2)| \le \frac{\ell_{\max}}{2\pi} |\th_1-\th_2|
+ 7 (\ell_{\max}^{3/4}/A^{1/4} + \ell_{\max}^{1/4}) \sqrt{L (t_1-t_2)}$$
whenever $|t_1-t_2| \le \min(2\sqrt{A\ell_{\min}},\ell_{\min}^{3/2})/(8L)$.
\end{proclaim}

\begin{demo}{Proof}
We need to compare the constant speed parametrization here with the horizontal
parametrization -- call it $c^*$ -- used in the maximum
distance bound \nmb!{3.5}. Under the horizontal
parametrization,
let the point $(t_1,\th_1)$ on $C_{t_1}$ correspond to $(t_2,\th_1^*)$ on
$C_{t_2}$, i.e.\ $c(t_2,\th_1^*) = c^*(t_2,\th_1)$. Let 
$C=(\ell_{\max}^{3/4}/A^{1/4} + \ell_{\max}^{1/4})$. Then we
know from \nmb!{3.5} that
$$ |c(t_1,\th_1) - c(t_2,\th_1^*)| \le 4C \sqrt{L(t_1-t_2)}.$$
To compare $\th_1$ and $\th_1^*$, we use the properties of the set $E$ in 
the proof of \nmb!{3.5} to estimate:
\begin{align*} \frac{(\th_1^* - \th_1) \ell_2}{2\pi} &=
\int_0^{\th_1} |c_\th^*(t_2,\ph)| d\ph - \frac{\th_1 \ell_2}{2\pi} \\
&\ge \Bigl(1-\frac{\sqrt{L(t_1-t_2)}}{(A\ell_1)^{1/4}} \Bigr) (\th_1 - \mu(E)) 
  \frac{\ell_1}{2\pi} -\frac{\th_1 \ell_2}{2\pi} \\
& \ge -2 \ell_1 \frac{\sqrt{L(t_1-t_2)}}{(A\ell_1)^{1/4}} - |\ell_1-\ell_2| 
  \mbox{ and similarly} \\
\frac{((2\pi - \th_1^*)-(2\pi-\th_1))\ell_2}{2\pi} &=
\int_{\th_1}^{2\pi} |c_\th^*(t_2,\ph)| d\ph - \frac{(2\pi-\th_1) \ell_2}{2\pi} \\
& \ge -2 \ell_1 \frac{\sqrt{L(t_1-t_2)}}{(A\ell_1)^{1/4}} - |\ell_1-\ell_2| 
\end{align*}
Combining these and using the Lipschitz property of length, we get:
\begin{align*} \frac{|\th_1^* - \th_1|\ell_2}{2\pi} &\le 
2C \sqrt{L(t_1-t_2)} + 
2|\sqrt{\ell_1}-\sqrt{\ell_2}| \sqrt{\ell_{\max}} \\
&\le 2C \sqrt{L(t_1-t_2)} + 
\sqrt{\ell_{\max}}\frac{L(t_1-t_2)}{\sqrt{A}} \le \frac{5}{2}
C \sqrt{L(t_1-t_2)} 
\end{align*}
Thus, finally:
\begin{align*}
|c(t_1,\th_1)-c(t_2,\th_2)| &\le |c(t_1,\th_1)-c(t_2,\th_1^*)| + 
\\&\qquad\qquad +
|c(t_2,\th_1^*)-c(t_2,\th_1)| + |c(t_2,\th_1)-c(t_2,\th_2)| 
\\&
\le 4C\sqrt{L(t_1-t_2)} + \frac{5}{2} C \sqrt{L(t_1-t_2)} 
+\frac{\ell_{\max}}{2\pi} |\th_1-\th_2|. \qed
\end{align*}
\end{demo}

\subsection*{\nmb.{3.9}}
One might also ask whether the maximum distance bound \nmb!{3.5} 
can be strengthened to assert that the 
1-jets of such curves $C$ must be close to the 1-jets of $C_0$. The
answer is NO, as is easily seen from looking a small wavelet-type perturbations 
of $C_0$. Specifically, calculate the length of the path:
$c(t,\th) = c_0(\th) + t\cdot a f(\th/a)\cdot i(c_0)_\th(\th), 0 \le t \le 1$ 
where $f(x)$ is an {\it arbitrary} $C^2$ function with compact support and
$a$ is very small. We claim the length of this path is $O(\sqrt{a})$, while
the 1-jet at the point $\th=0$ of the final curve of the path approaches
$(1+if'0)(c_0)_\th (0)$.

We sketch the proof, which is straightforward. Let $C_{a,t}$ be the curves
on this path. Then $\sup |c_t| = O(a),\, \sup |\kappa_{C_{a,t}}| = O(1/a),\,
A \le |c_\th| \le B$ for suitable $A,B>0$
and $\ell(\on{support}(c_t)) = O(a)$. Then the integral 
$\int_{S^1} (1+A\kappa_c^2) (c_t,ic_\th)^2 \frac{d\th}{|c_\th|}$ breaks up
into 2 pieces, the first being $O(a^2)$, the second being $O(1)$ and the
integral vanishing outside an interval of length $O(a)$. Thus the total
distance is $O(\sqrt{a})$.

\subsection*{\nmb.{3.10}. The $H^0$-distance on $B_i(S^1,\mathbb R^2)$ vanishes 
}

Let $c_0,c_1\in\on{Imm}(S^1,\mathbb R^2)$ be two immersions, 
and suppose that 
$t\mapsto (\th\mapsto c(t,\th))$ is a smooth curve in $\on{Imm}(S^1,\mathbb R^2)$
with $c(0,\cdot)=c_0$ and $c(1,\cdot)=c_1$. 

The arc-length for the $H^0$-metric of the curve 
$t\mapsto \pi(c(t,\cdot))$ in $B_i(S^1,\mathbb R^2)$ is given by \nmb!{3.2.7}
as
\begin{equation}
L^{\text{hor}}_{G^0}(c) 
=\int_0^1\Bigl(\int_{S^1}
\langle c_t,ic_\th \rangle^2 \frac{d\th}{|c_\th|} \Bigr)^{\tfrac12}dt
\AMStag{\nmb:{1}}
\end{equation}

\begin{proclaim}{Theorem} 
For $c_0,c_1\in\on{Imm}(S^1,\mathbb R^2)$ there exists always
a path $t\mapsto c(t,\cdot)$ with 
$c(0,\cdot)=c_0$ and $\pi(c(1,\cdot))=\pi(c_1)$ 
such that $L^{\text{hor}}_{G^0}(c)$ is arbitrarily small.
\end{proclaim}

Heuristically, the reason for this is that if the curve is made to zig-zag
wildly, say with teeth at an angle $\al$, then the length of the curve
goes up by a factor $1/\cos(\al)$ but the {\it normal} component of the
motion of the curve goes down by the factor $\cos(\al)$ -- and this
normal component is squared, hence it dominates.

\begin{demo}{Proof}
Take a path $c(t,\th)$ in $\on{Imm}(S^1,\mathbb R^2)$ from $c_0$ to $c_1$ and
make it horizontal using \nmb!{2.5} so that
that $\langle c_t,c_\th \rangle=0$; this forces a reparametrization on
$c_1$. 

\begin{figure} \begin{center}
\epsfig{width=3in,file=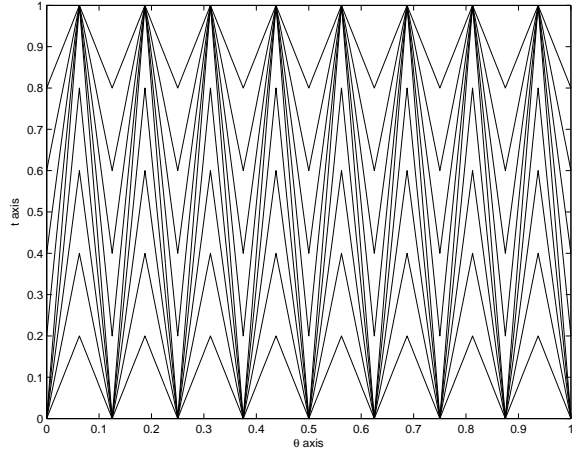}
\end{center} 
\caption{The reparametrization of a path of curves used to make its 
  length arbitrarily small.}
\end{figure}

Now let us view $c$ as a smooth mapping $c:[0,1]\x [0,1]\to \mathbb R^2$. 
We shall use the piecewise linear reparameterization
$(\ph(t,\th),\th)$ of the square shown above, which for $0\le t\le 1/2$ 
deforms the straight line into a zig-zag
of height $1$ and period $n/2$ connecting the two end-curves, and then 
removes the teeth for $1/2\le t\le 1$. In
detail: Let $\tilde c(t,\th)=c(\ph(t,\th),\th)$ where 
\begin{displaymath}
\ph(t,\th)=
\begin{cases}
  2t(2n\th-2k) &\text{ for }\; 0 \le t\le 1/2,\;
    \tfrac{2k}{2n}\le \th \le\tfrac{2k+1}{2n}
  \\
  2t(2k+2-2n\th) &\text{ for }\; 0 \le t\le 1/2,\;
    \tfrac{2k+1}{2n}\le \th \le\tfrac{2k+2}{2n}
  \\
  2t-1+2(1-t)(2n\th-2k) &\text{ for }\; 1/2 \le t\le 1,\;
    \tfrac{2k}{2n}\le \th \le\tfrac{2k+1}{2n}
  \\
  2t-1+2(1-t)(2k+2-2n\th) &\text{ for }\; 1/2 \le t\le 1,\;
    \tfrac{2k+1}{2n}\le \th \le\tfrac{2k+2}{2n}.
\end{cases}
\end{displaymath}
Then we get $\tilde c_\th = \ph_\th.c_t +c_\th$ and 
$\tilde c_t = \ph_t. c_t$ where
\begin{displaymath}
\ph_\th=\begin{cases} +4nt \\ -4nt \\+4n(1-t) \\ -4n(1-t) \end{cases},\qquad 
\ph_t  =\begin{cases} 4n\th-4k \\ 4k+4-4n\th \\2-4n\th+4k\\-(2-4n\th+4k)\end{cases}.
\end{displaymath}
Also, $\langle c_t,c_\th \rangle=0$ implies
$\langle \tilde c_t, i\tilde c_\th \rangle = \ph_t.|c_t|.|c_\th|$ and 
$|\tilde c_\th|=|c_\th|\sqrt{1+\ph_\th^2(|c_t|/|c_\th|)^2}$.
Thus 
\allowdisplaybreaks\begin{align*}
&L^{\text{hor}}(\tilde c)  
=\int_0^1\Bigl(\int_0^1
\langle \tilde c_t,i\tilde c_\th \rangle^2 
  \frac{d\th}{|\tilde c_\th|} \Bigr)^{\tfrac12}dt
=\int_0^1\Bigl(\int_0^1 
\frac{\ph_t^2 |c_t|^2 |c_\th|}
{\sqrt{1+\ph_\th^2(\frac{|c_t|}{|c_\th|})^2}}d\th\Bigr)^{\tfrac12}dt
=\\&
=\int_0^{\frac12}\Biggl(\sum_{k=0}^{n-1}\Bigl(
\int_{\tfrac{2k}{2n}}^{\tfrac{2k+1}{2n}}
  \frac{(4n\th-4k)^2 |c_t(\ph,\th)|^2 |c_\th(\ph,\th)|}
  {\sqrt{1+(4nt)^2(\frac{|c_t(\ph,\th)|}{|c_\th(\ph,\th)|})^2}}d\th
+\\&\qquad\qquad\qquad\qquad
+\int_{\tfrac{2k+1}{2n}}^{\tfrac{2k+2}{2n}}
  \frac{(4k+4-4n\th)^2 |c_t(\ph,\th)|^2 |c_\th(\ph,\th)|}
  {\sqrt{1+(4nt)^2(\frac{|c_t(\ph,\th)|}{|c_\th(\ph,\th)|})^2}}d\th
\Bigr)\Biggr)^{\tfrac12}dt
+\\&\quad
+\int_{\frac12}^1\Biggl(\sum_{k=0}^{n-1}\Bigl(
\int_{\tfrac{2k}{2n}}^{\tfrac{2k+1}{2n}}
  \frac{(2-4n\th+4k)^2 |c_t(\ph,\th)|^2 |c_\th(\ph,\th)|}
  {\sqrt{1+(4n)^2(1-t)^2(\frac{|c_t(\ph,\th)|}{|c_\th(\ph,\th)|})^2}}d\th
+\\&\qquad\qquad\qquad\qquad
+\int_{\tfrac{2k+1}{2n}}^{\tfrac{2k+2}{2n}}
  \frac{(2-4n\th+4k)^2 |c_t(\ph,\th)|^2 |c_\th(\ph,\th)|}
  {\sqrt{1+(4n)^2(1-t)^2(\frac{|c_t(\ph,\th)|}{|c_\th(\ph,\th)|})^2}}d\th
\Bigr)\Biggr)^{\tfrac12}dt
\end{align*}
The function $|c_\th(\ph,\th)|$ is uniformly bounded above and away from 0, 
and $|c_t(\ph,\th)|$ is uniformly bounded. Thus we may estimate
\begin{align*}
&\sum_{k=0}^{n-1}
\int_{\tfrac{2k}{2n}}^{\tfrac{2k+1}{2n}}
  \frac{(4n\th-4k)^2 |c_t(\ph,\th)|^2 |c_\th(\ph,\th)|}
  {\sqrt{1+(4nt)^2(\frac{|c_t(\ph,\th)|}{|c_\th(\ph,\th)|})^2}}d\th
\\&
\le O(1) \sum_{k=0}^{n-1}
\int_{0}^{\tfrac{1}{2n}}
  \frac{4n^2\th^2 |c_t(\ph(t,\frac{2k}{2n}+\th),\tfrac{2k}{2n}+\th)|^2 }
  {\sqrt{1+(4nt)^2|c_t(\ph(t,\frac{2k}{2n}+\th),\tfrac{2k}{2n}+\th)|^2}}d\th
\end{align*}
We estimate as follows.
Fix $\ep>0$. First we split of the integral $\int_{t=0}^\ep$
which is $O(\ep)$ uniformly in $n$; so for the rest we have
$t\ge \ep$. The last sum of integrals is now estimated as follows:
Consider first the set of all $\th$ such that 
$|c_t(\ph(t,\frac{2k}{2n}+\th),\tfrac{2k}{2n}+\th)|<\ep$ which is
a countable disjoint union of open intervals. There we get the estimate
$O(1).n. 4n^2. \ep^2 (\th^3/3)|_{\th=0}^{\th=1/2n}=O(\ep)$, uniformly in $n$.
On the complementary set of all $\th$ where 
$|c_t(\ph(t,\frac{2k}{2n}+\th),\tfrac{2k}{2n}+\th)|\ge\ep$
we use also $t\ge \ep$ and estimate by
$O(1).n.4n^2.\tfrac{1}{4n\ep^2}.(\th^3/3)|_{\th=0}^{\th=1/2n}=O
(\tfrac{1}{\ep^2 n})$.
The other sums of integrals can be estimated similarly, thus 
$L^{\text{hor}}(\tilde c)$ goes to 0 for $n\to \infty$. It is clear that
one can approximate $\ph$ by a smooth function without changing the
estimates essentially. 
\qed\end{demo}

\subsection*{\nmb.{3.11}. Non-smooth curves in the completion of $B_i$.}
We have seen in \nmb!{3.6} that the completion of $B_i$ in the metric $G^A$ lies in the
space of Lipschitz maps $c: S^1 \rightarrow \mathbb R^2$ mod monotone 
correspondences, that is, rectifiable Fr\'echet immersed curves. But how big
is it really? We cannot answer this, but we show, in this section, that certain
non-smooth curves are in the completion. To be precise, if $c$ is rectifiable,
then we can assume $c$ is parametrized at constant speed $|c_\th| \equiv L/2\pi$
where $L$ is the length of the curve. Therefore $c_\th = (L/2\pi)e^{i\al(\th)}$ 
for some measurable function $\al(\th)$ giving the orientation of the tangent 
line at almost every point. We will say that a rectifiable curve $c$ 
is {\it 1-BV} if the function $\al$ is of bounded variation. Note that this
means that the derivative of $\al$ exists as a finite signed measure, hence
the curvature of $c$ -- which is $(2\pi/L) \al'$ -- is also a finite signed 
measure.
In particular, there are a countable set of `vertices' on such a curve, points
where $\al$ has a discontinuity and the measure giving its curvature has an
atomic component. Note that $\al$ has left and right limits everywhere and
vertices can be assigned angles, namely $\al_+(\th)-\al_-(\th)$.

\begin{proclaim}{Theorem}
All 1-BV rectifiable curves are in the completion of $B_i$ with respect to the 
metric $G^A$.
\end{proclaim}

\begin{demo}{Proof}
This is proven using the following lemma:
\begin{proclaim}{Lemma}
Let $c(t,\th), \, 0 < t \le 1$ be an open path of smooth curves $c(t)$ and let
$\al(t,\th) = \arg(c_\th(\th,t))$. Assume that
\begin{enumerate}
\item the length of all curves $c(t)$ is bounded by $C_1$,
\item $|c_t| \le C_2$, for all $(t,\th)$,
\item For all $t$, the total variation in $\th$ of $\al(\th,t)$ is
bounded by $C_3$ and
\item the curvature of $c(t)$ satisfies $|\ka_{c(t)}(\th,t)| \le C_4/t$
for all $\th$.
\end{enumerate}
Then the length of this path is bounded by $C_2(\sqrt{C_1}+2\sqrt{AC_3C_4})$.
\end{proclaim}

To prove the lemma, let $s_t$ be arc length on $c(t)$ and estimate the integral:
\begin{align*}
\int_{c(t)} (1+A\ka(c(t))(t,\th)^2) 
  \langle c_t, \frac{ic_\th}{|c_\th|} \rangle^2 |c_\th| d\th 
  &\le C_2^2 \bigl( C_1 + A\int_{c(t)} \ka_{c(t)}^2 ds_t \bigr)
\\& 
= C_2^2 \bigl( C_1 + A\int_{c(t)} \ka_{c(t)} \frac{d\al}{ds_t} ds_t \bigr) 
\\& 
\le C_2^2 \bigl( C_1 + A \frac{C_4}{t} C_3\bigr).
\end{align*}
Taking the square root of both sides and integrating from 0 to 1, we get the
result.
\end{demo}

We apply this lemma to the simplest possible smoothing of a 1-BV rectifiable
curve $c_0$:
$$ c(t,\th) = \frac{1}{\sqrt{2\pi}t} \int_{\mathbb R} 
c_0(\th-\ph) e^{-\ph^2/2t^2} d\ph = \frac{1}{\sqrt{2\pi}t} \int_{\mathbb R} 
c_0(\ph) e^{-(\th-\ph)^2/2t^2} d\ph, 0 < t \le 1.$$
Note that $t$ is the standard deviation of the Gaussian, {\it not} the 
variance. We assume $c_0$ has a constant speed parametrization and 
$c'_0 = (L/2\pi)e^{i\al}$ as above, where $\al'$ is a finite signed measure. Thus:
\begin{align*}
c_\th &= \frac{L}{(2\pi)^{3/2}t} \int_{\mathbb R} e^{i\al(\th-\ph)-\ph^2/2t^2} d\ph 
\\
c_{\th \th} &= \frac{iL}{(2\pi)^{3/2}t} \int_{\mathbb R} e^{i\al(\ph) - 
(\th-\ph)^2/2t^2} \al'(d\ph)
\end{align*}
Moreover, using the second expression for the convolution and the heat 
equation for the Gaussian, we see that $c_t = tc_{\th \th}$. We now estimate:
\begin{align*}
|c_\th| &\le L/2\pi, \quad \mbox{hence length} (C_t) \le L\\
|c_{\th \th}| &\le \frac{L}{(2\pi)^{3/2}t} \int_{S^1} 
\sum_n e^{-(\th-\ph-nL)^2/2t^2} |\al'|(d\ph) \\
&\le \sup_x \bigl( \sum_n e^{-(x-nL)^2/2t^2}\bigr) 
  \frac{L \cdot \mbox{Var}(\arg(c'_0))}{(2\pi)^{3/2}t} = O(1/t), \\
\int_{S^1} |c_{\th \th}| d\th &\le \frac{L}{2\pi} 
  \Bigl(\int_{\mathbb R}\frac{1}{\sqrt{2\pi}t} e^{-\th^2/2t^2} d\th \Bigr) 
  \Bigl( \int_{S^1} |\al'(d\ph)|\Bigr) 
  = \frac{L}{2\pi} \mbox{Var}(\arg (c'_0)) \\
|c_t| &= t|c_{\th \th}| = O(1).
\end{align*}

To finish the proof, all we need to do is get a lower bound on $|c_\th|$.
However, $|c_\th|$ can be very small if the curve $c_0$ has corners with
small angles. In fact, $c_0$ can even double back on itself, giving a `corner'
with angle $\pi$. We need to treat this as a special case. When all the vertex 
angles of $c_0$ are less than $\pi$, we can get a lower bound for $|c_\th|$ as
follows. We start with the estimate:
\begin{align*}
|c_\th(\th)| &= \Bigl| \frac{1}{\sqrt{2\pi}t} \int_{\mathbb R}
e^{i\al(\th-\ph) - \ph^2/2t^2} d\ph \bigl/ e^{i\al(\th)} \Bigr| \\
&\ge \Bigl| \frac{1}{\sqrt{2\pi}t} \int_{\mathbb R} \cos(\al(\th-\ph)-\al(\th))
e^{-\ph^2/2t^2} d\ph \Bigr|
\end{align*}
We break up the integral over $\mathbb R$ into 3 intervals 
$(-\infty, \th-\de/2], [\th-\de/2, \th+\de/2], [\th+\de/2, +\infty)$ for a
suitable $\de$. If
$t$ is sufficiently small, the integral of the Gaussian over the first and 
third intervals goes uniformly to 0 and, on the middle interval, goes to 1.
Thus it suffices to estimate the $\cos$ in the middle interval. We use a
remark on BV functions:

\begin{proclaim}{Lemma}
For any BV function $f(x)$ and any $C>0$, there is a $\delta>0$ such that
on every interval $I$ of length less than $\delta$, either $f|_I$ has a
single jump of size $\ge C$ or $\max(f|_I) - \min(f|_I) \le C$. 
\end{proclaim}

In fact, let $C-\ep$ be the size of the largest 
jump in $f$ less than $C$
and break up the domain of $f$ into intervals $J_i$ on each of which the
variation of $f$ is less than $\epsilon/2$, big jumps being on their
boundaries. If $\delta$ is less than the minimum of the lengths of the $J_i$,
we get what we want.

Now let $\pi - \be$ be the largest vertex angle of the curve $c_0$. Using 
the last lemma, choose a $\de$ so that on every interval $I$ in the $\th$-line of
length less than $\de$, either $I$ contains a single vertex with exterior angle 
$\ge \be/3$ or $\max \al|_I - \min \al|_I \le \be/3$. Now if there is no
vertex in $[\th-\de/2, \th+\de/2]$, then $|\al(\th-\ph)-\al(\th)| \le \be/3$
on this interval and our lower bound is:
$$ |c_\th(\th)| \ge \cos(\be/3) - o(t).$$

On the other hand , if there is such a vertex, say at $\bar \th$, then 
$\al$ varies by at most $\be/3$ in $[\th-\de/2,\bar \th)$, jumps by at
most $\pi-\be$ at $\bar \th$ and then varies by at most $\be/3$ on
$(\bar \th, \th+\de/2]$. Assume $\th < \bar \th$ (the case $\th > \bar \th$
is similar). Then:
\begin{displaymath}
\cos(\al(\th-\ph)-\al(\th))\ge
\begin{cases}
\cos(\be/3), &\mbox{ if } \ph \in (\th-\bar \th, \th+\de/2] \\
\cos(\pi-\be+\be/3) = -\cos(2\be/3), &\mbox{ if } \ph \in 
[\th-\de/2,\th- \bar \th)
\end{cases}
\end{displaymath}
Thus:
$$ |c_\th(\th)| \ge \tfrac{1}{2} (\cos(\be/3)-\cos(2\be/3)) - o(t).$$
hence, if $t$ is sufficiently small, we get a uniform lower bound on $|c_\th|$.
Since $|\ka_{C_t}| \le |c_{\th \th}|/|c_\th|^2$, we get the required upper
bound both on $|\ka_{C_t}|$ and on the variation of $\al_{C_t}$, i.e.\ 
$\int_{S^1} |\ka_{C_t}|$ and all the requirements of the lemma are satisfied.

\begin{figure}[htb]
\begin{center}
\epsfig{width=3in,file=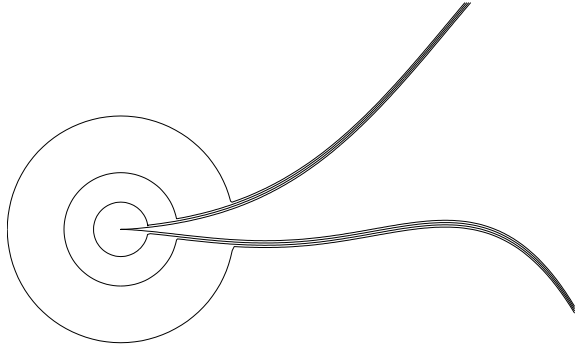}
\caption{Approximating $1-BV$ curves with zero angle vertices by curves
with positive angle vertices.}
\end{center}
\end{figure}

If $c_0$ has a vertex with angle $\pi$, we need to add an extra argument. $c_0$
certainly has at most a finite number of such vertices and we can construct a
new curve by drawing a circle of radius $t$ around each of these vertices and
letting $c_0^{(t)}$ be the curve which follows $c_0$ until it hits one of these
circles and then replaces the vertex with a circuit around the circle: see Figure 3. 
Each of the curves $c_0^{(t)}$ is in the completion of $B_i$ by the previous
argument and the path formed by the $c_0^{(t)}$'s also has finite length,
hence $c_0$ is in the completion. We omit the details which are straightforward.

\subsection*{\nmb.{3.12}. The energy of a path as `anisotropic area' 
of its graph in $\mathbb R^3$}

Consider a path $t\mapsto c(t,\cdot)$ in the manifold $\on{Imm}(S^1,\mathbb R^2)$.
It projects to a path $\pi \circ c$ in $B_i(S^1,\mathbb R^2)$ whose energy is
\begin{align}
E_{G^A}(\pi\o c) 
&= \tfrac12\int_a^b G^A_{\pi(c)}(T_c\pi.c_t,T_c\pi.c_t)\,dt
\notag\\& 
= \tfrac12\int_a^b G^A_{c}(c_t^\bot,c_t^\bot)\,dt
=\tfrac12\int_a^b\int_{S^1}(1+A\ka(c)^2) 
  \langle c_t^\bot,c_t^\bot\rangle |c_\th|\,d\th\,dt
\notag\\&
=\tfrac12\int_a^b\int_{S^1}(1+A\ka(c)^2)
\left\langle \frac{\langle c_t,ic_\th \rangle}{|c_\th|^2}ic_\th, 
\frac{\langle c_t,ic_\th \rangle}{|c_\th|^2}ic_\th \right\rangle
|c_\th|\,d\th\,dt
\notag\\&
=\tfrac12\int_a^b\int_{S^1} (1+A\ka(c)^2)
\langle c_t,ic_\th \rangle^2 \frac{d\th}{|c_\th|}\,d\th dt
\AMStag{\nmb:{1}}
\end{align}
If the path $c$ is horizontal, i.e., it
satisfies $\langle c_t,c_\th \rangle=0$. Then $\langle  c_t,ic_\th\rangle=
|c_t|.|c_\th|$ and we have
\begin{equation}
E^{\text{hor}}_{G^A}(c) 
=\tfrac12\int_a^b\int_{S^1}(1+A\ka(c)^2)|c_t|^2|c_\th|\,d\th\, dt,
\qquad \langle c_t,c_\th \rangle=0
\AMStag{\nmb:{2}}
\end{equation}
which is just the usual energy of $c$.

Let $c(t,\th)=(x(t,\th),y(t,\th))$ be still horizontal and consider the graph 
$$\Ph(t,\th)=(t,x(t,\th),y(t,\th)) \in \mathbb R^3.$$
We also have $|x_ty_\th-x_\th
y_t|=|\det(c_t,c_\th)|=|c_t|.|c_\th|$ and for the vector product
$\Ph_t\x \Ph_\th = (x_ty_\th-x_\th y_t, -y_\th,x_\th)$, so we get 
$$|\Ph_t\x \Ph_\th|^2 = (x_ty_\th-x_\th y_t)^2 + y_\th^2 +x_\th^2
  = (x_\th^2 + y_\th^2)(x_t^2+y_t^2 + 1) = |c_\th|^2(|c_t|^2+1).$$
We express now $E^{\text{hor}}(c)$ as an integral over the
immersed surface
$S\subset \mathbb R^3$ parameterized by $\Ph$ in terms of the surface area
$d\mu_S = |\Ph_t\x \Ph_\th|\,d\th\,dt$ as follows: 
\begin{align*}
E^{\text{hor}}_{G^A}(c) &=\tfrac12\int_a^b\int_{S^1}(1+A\ka(c)^2)
  \frac{|c_t|^2|c_\th|}{|\Ph_t\x\Ph_\th|}|\Ph_t\x\Ph_\th|\,d\th\,dt
\\&
=\tfrac12\int_{[a,b]\x S^1}(1+A\ka(c)^2)
\frac{|c_t|^2}{\sqrt{|c_t|^2+1}}\,d\mu_S
\end{align*}
Next we want to express the integrand as as a function $\ga$ of the 
unit normal $n_S=(\Ph_t\x \Ph_\th)/|\Ph_t\x \Ph_\th|$. Let $e_0=(1,0,0)$,
then the absolute value of the $t$-component $n_S^0$ of the unit normal $n_S$ is
\begin{displaymath}
|n_S^0|:=|\langle e_0,n_S \rangle| 
= \frac{|c_t|}{\sqrt{|c_t|^2+1}},
\quad\text{ and }\quad
\frac{|c_t|^2}{\sqrt{|c_t|^2+1}}=\frac{|n_S^0|^2}{\sqrt{1-|n_S^0|^2}}.
\end{displaymath}
Thus for horizontal $c$ (i.e., with $c_t\bot c_\th$) we have 

\begin{proclaim}{Horizontal energy as anisotropic area}
\begin{equation}
\boxed{\quad
E^{\text{hor}}_{G^A}(c) 
=\tfrac12\int_{[a,b]\x S^1}(1+A\ka(c)^2)
  \frac{|n_S^0|^2}{\sqrt{1-|n_S^0|^2}}\,d\mu_S
  \quad}
\AMStag{\nmb:{3}}
\end{equation}
\end{proclaim}

Here the final expression is only in terms of the surface $S$ and does not
depend on the curve $c$ being horizontal. 
This anisotropic area functional has to be minimized in order to prove that 
geodesics exists between arbitrary curves (of the same degree) in
$B_{i}(S^1,\mathbb R^2)$.
Thus we are led to the 

\begin{proclaim}{Question} For immersions $c_0,c_1:S^1\to \mathbb R^2$ does
there exist an immersed surface 
$S=(\on{ins}_{[0,1]}, c):[0,1]\x S^1\to \mathbb R\x \mathbb R^2$
such that the functional \thetag{\nmb|{3}} is critical at $S$?
\end{proclaim}

A first step is:

\begin{proclaim}{Bounding the area} 
For any path $[a,b]\ni t\mapsto c(t,\quad)$ the area of the graph surface $S=S(c)$
is bounded as follows:
\begin{equation}
\text{Area}(S) = \int_{[a,b]\x S^1} d\mu_S \le
2 E^{\text{hor}}_{G^A}(c)+\max_{t}\ell(c(t,\quad))(b-a)
\AMStag{\nmb:{4}}
\end{equation}
\end{proclaim}

\begin{demo}{Proof}
Writing the unit normal $n_S=(n_S^0,n_S^1,n_S^2)\in S^2$ according to the
coordinates $(t,x,y)$ we have 
$$
|n_S^1|+|n_S^2|+\frac{|n_S^0|^2}{\sqrt{1-|n_S^0|^2}}
\ge |n_S^1|^2+|n_S^2|^2+|n_S^0|^2 = 1
$$
Since $|n_S^1|d\mu_S$ is the area element of the projection of $S$ onto the
$(t,y)$-plane we have
\begin{align*}
\text{Area}(S)
&= 
\int_{[a,b]\x S^1}\!\! d\mu_S 
\le \int_{[a,b]\x S^1}\!\! (1+A\ka(c)^2)
  \Bigl(|n_S^1|+|n_S^2|+\frac{|n_S^0|^2}{\sqrt{1-|n_S^0|^2}}\Bigr)d\mu_S 
\\&\le
  2 E^{\text{hor}}_{G^A}(c)+\max_{t}\ell(c(t,\quad))(b-a).
\qed
\end{align*}
\end{demo}

\section*{\totoc\nmb0{4}. Geodesic equations and sectional curvatures 
}

\subsection*{\nmb.{4.1}. Geodesics on $\on{Imm}(S^1,\mathbb R^2)$ }
The energy of a curve $t\mapsto c(t,\cdot)$ in the space 
$\on{Imm}(S^1,\mathbb R^2)$ is 
\begin{displaymath}
E_{G^A}(c) = \tfrac12\int_a^b\int_{S^1} 
     (1+A\ka_c^2)\langle c_t,c_t\rangle |c_\th|\,d\th\,dt. 
\end{displaymath}
By calculating its first variation, we get the equation for a geodesic:
\begin{proclaim}{Geodesic Equation}
\begin{equation}\boxed{\quad
\Bigl( (1+A\ka^2)|c_\th|\cdot c_t \Bigr)_t = \Bigl( \frac{-1+A\ka^2}{2}\cdot
\frac{|c_t|^2}{|c_\th|}\cdot c_\th 
+ A \frac{(\ka |c_t|^2)_\th}{|c_\th|^2}\cdot i c_\th \Bigr)_\th.
\quad}
\AMStag{\nmb:{1}}
\end{equation}
\end{proclaim}

\begin{demo}{Proof}
 From \nmb!{2.2} we have
\begin{align*}
\ka(c)_s &= \frac{\langle ic_{s\th},c_{\th\th} \rangle}{|c_\th|^3}
     + \frac{\langle ic_\th,c_{s\th\th} \rangle}{|c_\th|^3}
     -3\ka\frac{\langle c_{s\th},c_\th \rangle}{|c_\th|^2}.
\end{align*}
and
\begin{align*}
c_{\th \th} &= \frac{\langle c_{\th \th}, c_\th \rangle}{|c_\th|^2} c_\th
+ \frac{\langle c_{\th \th}, ic_\th \rangle}{|c_\th|^2} ic_\th
\\ &= \frac{|c_\th|_\th}{|c_\th|} c_\th + \ka(c) |c_\th| ic_\th.
\end{align*}
Now we compute
\allowdisplaybreaks\begin{align*}
&\p_s|_0 E(c) = \tfrac12\p_s|_0 \int_a^b\int_{S^1} 
     (1+A\ka^2)\langle c_t,c_t \rangle  |c_\th|\,d\th\,dt 
\\&= \int_a^b \int_{S^1}\Bigl(
A\ka\ka_s|c_\th\|c_t|^2 
+ (1+A\ka^2)\langle c_{st},c_t \rangle |c_\th|
+\frac{1+A\ka^2}{2} |c_t|^2 
\frac{\langle c_{s\th},c_\th \rangle}{|c_\th|}
\Bigr)\,d\th dt
\\&= \int_a^b \int_{S^1}\Bigl(
 A\ka\langle ic_{s\th},c_{\th\th} \rangle\frac{|c_t|^2}{|c_\th|^2}
 +A\ka\langle ic_\th,c_{s\th\th} \rangle\frac{|c_t|^2}{|c_\th|^2}
 -3A\ka^2\langle c_{s\th},c_\th \rangle\frac{|c_t|^2}{|c_\th|} 
\\&\qquad\qquad\qquad
 -\Bigl\langle c_{s},
  \bigl( (1+A\ka^2) |c_\th| c_t \bigr)_t  
  +\bigl(\frac{1+A\ka^2}{2} \frac{|c_t|^2}{|c_\th|}c_\th\bigr)_\th 
 \Bigr\rangle
\Bigr)\,d\th dt
\\&= \int_a^b \int_{S^1}\Bigl(
 \Bigl\langle c_{s},
     A\Bigl(\ka\frac{|c_t|^2}{|c_\th|^2} ic_{\th\th}\Bigr)_\th 
  \Bigr\rangle
 +\Bigl\langle c_{s},
     A\Bigl(\ka \frac{|c_t|^2}{|c_\th|^2}ic_\th\Bigl)_{\th\th} 
  \Bigr\rangle
 +\Bigl\langle c_{s},
     3A\Bigl(\ka^2\frac{|c_t|^2}{|c_\th|} c_\th \Bigr)_\th
   \Bigr\rangle
\\&\qquad\qquad\qquad
 -\Bigl\langle c_{s},
    \bigl( (1+A\ka^2) |c_\th| c_t \bigr)_t  
  +\bigl(\frac{1+A\ka^2}{2} \frac{|c_t|^2}{|c_\th|}c_\th\bigr)_\th 
  \Bigr\rangle
\Bigr)\,d\th dt
\\&= \int_a^b \int_{S^1}\Bigl\langle c_{s},
 -\bigl((1+A\ka^2)|c_\th|c_t \bigr)_t + F_\th
\Bigr\rangle\,d\th dt
\end{align*}
%\newpage
where 
\begin{align*} F &= A\ka \frac{|c_t|^2}{|c_\th|^2}ic_{\th \th} 
+A(\ka |c_t|^2)_\th \frac{ic_\th}{|c_\th|^2}
-2A\ka |c_t|^2\frac{|c_\th|_\th ic_\th}{|c_\th|^3}
+A\ka \frac{|c_t|^2}{|c_\th|^2}ic_{\th \th} +
\\ &\qquad\qquad +3A\ka^2 \frac{|c_t|^2}{|c_\th|} c_\th 
  - \frac{1+A\ka^2}{2} \frac{|c_t|^2}{|c_\th|}c_\th 
\end{align*}
Substituting the expression for $c_{\th \th}$ and simplifying, this reduces to
$$ 
F = \frac{-1+A\ka^2}{2} \frac{|c_t|^2}{|c_\th|} c_\th
+A(\ka |c_t|^2)_\th \frac{i c_\th}{|c_\th|^2}
$$
which gives the required formula for geodesics.
\end{demo}

Putting $A=0$ in \nmb!{4.1.1} we get the  
geodesic equation for the $H^0$-metric on $\on{Imm}(S^1,\mathbb R^2)$
\begin{equation}
\bigl(|c_\th| c_t\bigr)_t 
= - \frac{1}{2}\bigl( \frac{|c_t|^2 c_\th}{|c_\th|} \bigr)_\th 
\AMStag{\nmb:{2}} 
\end{equation}

\subsection*{\nmb.{4.2} Geodesics on $B_i(S^1,\mathbb R^2)$}
We may also restrict to geodesics which are perpendicular to the orbits of
$\on{Diff}(S^1)$, i.e.\ $\langle c_t, c_\th \rangle \equiv 0$, obtaining the 
geodesics in the quotient space $B_i(S^1,\mathbb R^2)$. To write this in the
simplest way, we introduce the `velocity' $a$ by setting $c_t = ia c_\th/|c_\th|$
(so that $|c_t|^2 = a^2$). When we substitute this into the above geodesic 
equation, the equation splits into a multiple of $c_\th$ and a multiple of
$ic_\th$. The former vanishes identically and the latter gives:
\begin{align*} \bigl( (1+A\ka^2) |c_\th| a \bigr)_t \frac{ic_\th}{|c_\th|} &=
\frac{-1+A\ka^2}{2} a^2 \bigl(\frac{c_\th}{|c_\th|}\bigr)_\th 
+A\bigl(\frac{(\ka a^2)_\th}{|c_\th|}\bigr)_\th \frac{ic_\th}{|c_\th|}, 
\qquad \text{or}\\
\bigl( (1+A\ka^2)|c_\th| a \bigr)_t &= 
\frac{-1+A\ka^2}{2} \ka |c_\th| a^2
+ A\bigl(\frac{(\ka a^2)_\th}{|c_\th|}\bigr)_\th.
\end{align*}
If we use derivatives with respect to arclength instead of $\th$ and write these 
with the subscript $s$, so that $f_s=f_\th/|c_\th|$, this simplifies. We need:
$$ |c_\th|_t = \frac{\langle c_\th, c_{t \th} \rangle}{|c_\th|} =
-\frac{\langle c_{\th \th}, c_t \rangle}{|c_\th|} =
-a \frac{\langle c_{\th \th}, ic_\th \rangle}{|c_\th|^2} =
-a \ka |c_\th|$$
which gives us a simple form for the equation for geodesics on 
$B_i(S^1,\mathbb R^2)$:
\begin{equation}\AMStag{\nmb:{1}}
\bigl( (1+A\ka^2) a\bigr)_t = \frac{1+3A\ka^2}{2} \ka a^2 + A (\ka a^2)_{ss}.
\end{equation}
Finally, we may expand the $t$-derivatives on the left hand side, 
using the formula $\ka_t = a\ka^2 + a_{ss}$ noted in \nmb!{2.2.7}; we also
collect all constraint equations that we chose along the way:
\begin{equation}
\boxed{\quad\begin{aligned}
0 &= \langle c_t, c_s \rangle,\quad c_t=aic_s,\quad \ka = \langle c_{ss},ic_s \rangle
\\
a_t &= \frac{\tfrac12 \ka a^2 + 
A\bigl( a^2(\ka_{ss} - \tfrac12 \ka^3)
+ 4\ka_s a a_s + 2\ka a_s^2\bigr)}{1+A\ka^2}.
\end{aligned}\quad}
\AMStag{\nmb:{2}}
\end{equation}
Handle this with care: Going to unit speed parametrization (so that $f_s$
is really a holonomic partial derivative) destroys the first constraint
`horizontality'. This should be seen as a gauge fixing. 

\subsection*{\nmb.{4.3}. Geodesics on $B_i(S^1,\mathbb R^2)$ for
$A=0$}
Let us now set $A=0$.
We keep looking at horizontal geodesics, so that 
$\langle c_t,c_\th \rangle=0$ and $c_t=iac_\th/|c_\th|$ for
$a\in C^\infty(S^1)$. We use the functions $a$, $s=|c_\th|$,
and $\ka$. We use equations from \nmb!{4.2} but we do not use
the anholonomic derivative:
\begin{align*}
s_t=-a\ka s,\qquad a_t=\tfrac12\ka a^2,\qquad
\ka_t = a\ka^2 + \frac1s\left(\frac{a_\th}{s} \right)_\th
= a\ka^2 + \frac{a_{\th\th}}{s^2} - \frac{a_\th s_\th}{s^3}. 
\AMStag{\nmb:{1}}
\end{align*}
We may assume that $s|_{t=0}$ is constant. Let $v(\th)=a(0,\th)$ be
the initial value for $a$. Then from equations \thetag{\nmb|{1}}
we get 
$$
\frac{s_t}{s} = - a\ka = - 2\frac{a_t}{a}\quad \implies\quad 
\log(sa^2)_t=0 
$$
so that $sa^2$ is constant in $t$, 
\begin{equation}
s(t,\th)a(t,\th)^2=s(0,\th)a(0,\th)^2=v(\th)^2,
\AMStag{\nmb:{2}}
\end{equation}
a smooth family of conserved quantities along the geodesic.
This leads to the substitutions
$$
s=\frac{v^2}{a^2},\qquad \ka = 2\frac{a_t}{a^2}
$$
which transform the last equation \thetag{\nmb|{1}} to
\begin{equation}
a_{tt} - 4\frac{a_t^2}{a} - \frac{a^6a_{\th\th}}{2v^4} +
\frac{a^6a_\th v_\th}{v^5} - \frac{a^5a_\th^2}{v^4} =0,\quad
a(0,\th)=v(\th),
\AMStag{\nmb:{3}}
\end{equation}
a nonlinear hyperbolic second order equation.
Note that \thetag{\nmb|{2}} implies that wherever $v=0$ then
also $a=0$ for all $t$. 
For that reason, let us transform equation \thetag{\nmb|{3}} into
a less singular form by substituting 
$a=vb$. Note that $b=1/\sqrt{s}$. The outcome is 
\begin{equation}
(b^{-3})_{tt} = -\frac{v^2}{2}(b^3)_{\th\th} - 2vv_\th(b^3)_\th -
\frac{3vv_{\th\th}}{2}b^3, \qquad b(0,\th)=1.
\AMStag{\nmb:{4}}
\end{equation}

\subsection*{\nmb.{4.4}. 
The induced metric on $B_{i,f}(S^1,\mathbb R^2)$ in a chart }
We also want to compute the curvature of $B_i(S^1,\mathbb R^2)$ in this metric.
For this, we need second derivatives and the most convenient way to calculate 
these seems to be to use a local chart. Consider the smooth principal bundle 
$\pi:\on{Imm}_{f}(S^1,\mathbb R^2)\to B_{i,f}(S^1,\mathbb R^2)$ with structure 
group $\on{Diff}(S^1)$ described in \nmb!{2.4.3}. We shall describe 
the metric in the following chart near $C\in B_{i,f}(S^1,\mathbb R^2)$:
Let $c\in\on{Imm}_{f}(S^1,\mathbb R^2)$ be parametrized by arclength with 
$\pi(c)=C$ of length $L$, with unit normal $n_c$. 
We assume that the parameter $\th$ runs in the scaled circle $S^1_L$ below.
As in the proof of
\nmb!{2.4.3} we consider the mapping
\begin{align*}
&\ps: C^\infty(S^1_L,(-\ep,\ep)) \to \on{Imm}_{f}(S^1_L,\mathbb R^2),
  \qquad \mathcal Q(c):= \ps(C^\infty(S^1_L,(-\ep,\ep)))\\
&\ps(f)(\th) = c(\th) + f(\th)n_c(\th) = c(\th) + f(\th)ic'(\th),\\
&\pi\o\ps :C^\infty(S^1_L,(-\ep,\ep)) \to B_{i,f}(S^1,\mathbb R^2),
\end{align*}
where $\ep$ is so small that $\ps(f)$ is an embedding for each $f$.
By \nmb!{2.4.3} the mapping $(\pi\o\ps)\i$ is a smooth chart 
on $B_{i,f}(S^1_L,\mathbb R^2)$. 
Note that: 
\allowdisplaybreaks\begin{align*}
\ps(f)' &= c' + f'ic' + fic'' = (1-f\ka_c)c' + f'ic' \\
\ps(f)'' &= c'' + f''ic' + 2f'ic'' + fic'''
     = -(2f'\ka_c+f\ka_c')c'+ (\ka_c+f''-f\ka_c^2)ic'  \\
n_{\ps(f)} &= 
     \frac1{\sqrt{(1-f\ka_c)^2+{f'}^2}}\Bigl((1-f\ka_c)ic'-f'c'\Bigr),\\
T_f\ps.h &= h.ic'\quad \in C^\infty(S^1,\mathbb R^2) = 
     T_{\ps(f)}\on{Imm}_{f}(S^1_L,\mathbb R^2)\\
%&=\frac{h(1-f\ka_c)}{\sqrt{(1-f\ka_c)^2+{f'}^2}}
%      \frac{(1-f\ka_c)ic'-f'c'}{\sqrt{(1-f\ka_c)^2+{f'}^2}}\\
%&\qquad\qquad+\frac{hf'}{(1-f\ka_c)^2+{f'}^2}\Bigl((1-f\ka_c)c'+f'ic'\Bigr)\\
&=\frac{h(1-f\ka_c)}{\sqrt{(1-f\ka_c)^2+{f'}^2}}n_{\ps(f)}
     +\frac{hf'}{(1-f\ka_c)^2+{f'}^2}\ps(f)',\\
(T_f\ps.h)^\bot &=\frac{h(1-f\ka_c)}{\sqrt{(1-f\ka_c)^2+{f'}^2}}n_{\ps(f)}
     \quad\in \mathcal N_{\ps(f)},\\
\ka_{\ps(f)}&= \frac{1}{((1-f\ka_c)^2+{f'}^2)^{3/2}}
     \langle i\ps(f)',\ps(f)'' \rangle \\
%&= \frac{\langle(1-f\ka_c)ic' - f'c',-(2f'\ka_c+f\ka_c')c'+ 
%     (\ka_c+f''-f\ka_c^2)ic'\rangle}
%     {((1-f\ka_c)^2+{f'}^2)^{3/2}} \\
%&= \frac{(1-f\ka_c)(\ka_c+f''-f\ka_c^2) + f'(2f'\ka_c+f\ka_c')}
%     {((1-f\ka_c)^2+{f'}^2)^{3/2}} \\ 
%&= \frac{\ka_c(1-f\ka_c)^2 +f''(1-f\ka_c) + 2{f'}^2\ka_c+ff'\ka_c'}
%     {((1-f\ka_c)^2+{f'}^2)^{3/2}} \\ 
&= \frac{\ka_c+f''-2f\ka_c^2-ff''\ka_c+f^2\ka_c^3 + 
     2{f'}^2\ka_c+ff'\ka_c'}{((1-f\ka_c)^2+{f'}^2)^{3/2}}
\end{align*}
Let $G^A$ denote also the induced metric on 
$B_{i,f}(S^1_L,\mathbb R^2)$. Since $\pi$ is a Riemannian submersion, 
$T_{\ps(f)}\pi:(\mathcal N_{\ps(f)},G^A_{\ps(f)}) 
     \to (B_{i,f}(S^1_L,\mathbb R^2), G^A_{\pi(\ps(f))})$
is an isometry. 
Then we compute for 
$f\in C^\infty(S^1_L,(-\ep,\ep))$ and $h,k\in C^\infty(S^1_L,\mathbb R)$
\allowdisplaybreaks\begin{align*}
&((\pi\o\ps)^*G^A)_f(h,k) = 
G^A_{\pi(\ps(f))}\Bigl(T_f(\pi\o\ps)h,T_f(\pi\o\ps)k\Bigr)\\
&=G^A_{\ps(f)}\Bigl((T_f\ps.h)^\bot,(T_f\ps.k)^\bot\Bigr)\\
&=\int_{S^1_L}(1+A\ka_{\ps(f)}^2)
     \Bigl\langle(T_f\ps.h)^\bot,(T_f\ps.k)^\bot\Bigr\rangle 
     |\ps(f)'|\;d\th \\
&=\int_{S^1_L}(1+A\ka_{\ps(f)}^2)
     \frac{hk(1-f\ka_c)^2}{\sqrt{(1-f\ka_c)^2+{f'}^2}}\;d\th\\
%&=\int_{S^1_L}\left(1+A\left(
%     \frac{\ka_c+f''-2f\ka_c^2-ff''\ka_c+f^2\ka_c^3 + 
%     2{f'}^2\ka_c+ff'\ka_c'}{((1-f\ka_c)^2+{f'}^2)^{3/2}}
%     \right)^2\right)
%     \frac{hk(1-f\ka_c)^2}{\sqrt{(1-f\ka_c)^2+{f'}^2}}\;d\th \\
%&=\int_{S^1_L}\left(1+A\frac{\left(
%     \ka_c(1-f\ka_c)^2 +f''(1-f\ka_c) + 2{f'}^2\ka_c+ff'\ka_c'
%     \right)^2}{((1-f\ka_c)^2+{f'}^2)^{3}}
%     \right)
%     \frac{hk(1-f\ka_c)^2}{\sqrt{(1-f\ka_c)^2+{f'}^2}}\;d\th \\
%&=\int_{S^1_L}hk\,\frac{\left(((1-f\ka_c)^2+{f'}^2)^{3}+A\left(
%     \ka_c(1-f\ka_c)^2 +f''(1-f\ka_c) + 2{f'}^2\ka_c+ff'\ka_c'
%     \right)^2\right)(1-f\ka_c)^2}
%     {((1-f\ka_c)^2+{f'}^2)^{7/2}}\;d\th
\end{align*}
This is the expression from which we have to compute the geodesic equation 
in the chart on $B_{i,f}(S^1_L,\mathbb R^2)$.

\subsection*{\nmb.{4.5}. Computing the Christoffel symbols in 
$B_{i,f}(S^1_L,\mathbb R^2)$ at $C=\pi(c)$ } 
We have to compute second derivatives in $f$ of the expression of the 
metric in \nmb!{4.2}. For that we expand the two main contributing 
expressions in $f$ to order 2, where we put $\ka=\ka_c$.
\begin{align*}
&\ka_{\ps(f)} = 
\\&
=(1-2f\ka+f^2\ka^2+{f'}^2)^{-3/2}
(\ka+f''-2f\ka^2-ff''\ka+f^2\ka^3+2{f'}^2\ka+ff'\ka')
\\
&= \ka+(f''+f\ka^2)+(f^2\ka^3+\tfrac12{f'}^2\ka+ff'\ka'+2ff''\ka)
     +O(f^3)
\\
&(1-f\ka)^2(1-2f\ka+f^2\ka^2+{f'}^2)^{-1/2} 
  =1-f\ka-\tfrac12{f'}^2+O(f^3)
\end{align*}Thus
\begin{multline*}
(1+A\ka_{\ps(f)}^2)
     \frac{(1-f\ka_c)^2}{\sqrt{(1-f\ka_c)^2+{f'}^2}}  
= 1+A\ka^2+2Af''\ka+ Af\ka^3 -f\ka -
\\ 
-\tfrac12{f'}^2 +Af^2\ka^4
+A\tfrac12{f'}^2\ka^2 +2Aff'\ka\ka' +A{f''}^2 +4Aff''\ka^2 
\end{multline*}
and finally 
\begin{multline}
G^A_f(h,k) = ((\pi\o\ps)^*G^A)_f(h,k) =
\\
=\int_{S^1_L}hk\,\Bigl(
(1+A\ka^2)+(2Af''\ka+Af\ka^3 -f\ka)+
-\tfrac12{f'}^2 
\\
+A(4ff''\ka^2+f^2\ka^4
+\tfrac12{f'}^2\ka^2 +2ff'\ka\ka' +{f''}^2) + O(f^3) 
\Bigr)\;d\th.
\AMStag{\nmb:{1}}\end{multline}
We differentiate the metric
\begin{multline*}
dG^A(f)(l)(h,k) 
  =\int_{S^1_L}hk\,\Bigl(
  2Al''\ka+(A\ka^3 -\ka)l +4Alf''\ka^2+4Afl''\ka^2 + 
\\
+2Afl\ka^4 +(A\ka^2-1)f'l' +2Alf'\ka\ka'+2Afl'\ka\ka' +2Af''l'' + O(f^2) 
\Bigr)\;d\th
\end{multline*}and compute the Christoffel symbol
\allowdisplaybreaks\begin{align*}
&-2G^A_f(\Ga_f(h,k),l) = -dG^A(f)(l)(h,k) +dG^A(f)(h)(k,l) +dG^A(f)(k)(l,h) 
\\&
=\int_{S^1_L}l\Bigl(
(A\ka^3 -\ka +2A\ka\ka'f' +4A\ka^2f'' +2A\ka^4f)k h 
\\&\qquad\qquad
+(2A\ka +4A\ka^2f +2Af'')(h''k+hk'') 
\\&\qquad\qquad
+(A\ka^2f'-f' +2A\ka\ka'f)(h'k+hk') 
+O(f^2)\Bigr)\;d\th
\\&\quad
-\int_{S^1_L}\Bigl(
l'(A\ka^2f'hk-f'hk+2A\ka\ka'fhk)
\\&\qquad\qquad
+l''(2A\ka hk +4A\ka^2fhk +2Af''hk) 
+O(f^2)\Bigr)\;d\th
\\&
=\int_{S^1_L}l\Bigl(
(A\ka^3 -\ka -2A\ka'')hk -4A\ka'(h'k+hk') -4A\ka h'k'
\\&\qquad\qquad
+(-2Af^{(4)}\!-f''+2A\ka^4f-6A\ka'{}^2f-6A\ka\ka''f-10A\ka\ka'f'+A\ka^2f'')hk
\\&\qquad\qquad
-(2f'+4Af'''+12A\ka\ka'f+6A\ka^2f')(h'k+hk')
\\&\qquad\qquad
-2(4A\ka^2f+2Af'')h'k'
+O(f^2)\Bigr)\;d\th.
\end{align*}
Thus 
\begin{align*}
&G^A_f(\Ga_f(h,k),l) = 
\\&
=\int_{S^1_L}l\Bigl(
(\tfrac12\ka -\tfrac12A\ka^3 +A\ka'')hk +2A\ka'(h'k+hk') +2A\ka h'k'
\\&\qquad\qquad
-(-Af^{(4)}-\tfrac12f''+A\ka^4f-3A\ka'{}^2f-3A\ka\ka''f-5A\ka\ka'f'
  +\tfrac12A\ka^2f'')hk
\\&\qquad\qquad
+(f'+2Af'''+6A\ka\ka'f+3A\ka^2f')(h'k+hk')
\\&\qquad\qquad
+(4A\ka^2f+2Af'')h'k'
+O(f^2)\Bigr)\;d\th.
\end{align*}
At the center of the chart, for $f=0$, we get
\begin{align*}
&G^A_0(\Ga_0(h,k),l) = 
\\&
=\int_{S^1_L}l\Bigl(
(\tfrac12\ka-\tfrac12A\ka^3 +A\ka'')hk +2A\ka'(h'k+hk') +2A\ka h'k'
\Bigr)\;d\th
\\&
=\int_{S^1_L}l\Bigl(\frac{(\tfrac12\ka-\tfrac12A\ka^3 +A\ka'')hk 
  +2A\ka'(h'k+hk') +2A\ka h'k'}{1+A\ka^2}\Bigr)(1+A\ka^2)\;d\th
\\&
=G^A_0\Bigl(\frac{(\tfrac12\ka-\tfrac12A\ka^3 +A\ka'')hk 
  +2A\ka'(h'k+hk') +2A\ka h'k'}{1+A\ka^2},l\Bigr)
\end{align*}
so that 
\begin{equation}
\Ga_0(h,k) = \frac{(\tfrac12\ka-\tfrac12A\ka^3 +A\ka'')hk 
  +2A\ka'(h'k+hk') +2A\ka h'k'}{1+A\ka^2}.
\AMStag{\nmb:{2}}
\end{equation}
Letting $h=k=f_t$, this leads to the geodesic equation, valid at $f=0$:
$$ f_{tt} = \frac{(\tfrac12\ka-\tfrac12A\ka^3 +A\ka'')f_t^2 
  +4A\ka'f_t f'_t +2A\ka (f'_t)^2}{1+A\ka^2}.$$
If we substitute $a$ for $f_t$ and $a_t$ for $f_{tt}$, this is the same as
the previous geodesic equation derived in \nmb!{4.2} by variational
methods. There is a subtle point here, however: why is it ok to identify
the second derivatives $a_t$ and $f_{tt}$ with each other? To check this
let $c(\th)+(ta_1(\th)+\frac{t^2}{2} a_2(\th))ic'(\th)$ be a 2-jet in our chart.
Then if we reparametrize the nearby curves by substituting 
$\th - \frac{t^2}{2}a_1 a'_1$ for $\th$, letting 
\begin{align*}
c(t,\th) &= c(\th - \tfrac{t^2}2 a_1 a'_1) 
  + (t a_1(\th - \tfrac{t^2}2 a_1 a'_1) 
  +\frac{t^2}{2}a_2(\th - \tfrac{t^2}2 a_1 a'_1))
  ic(\th - \tfrac{t^2}2 a_1 a'_1)'\\
&\equiv c(\th)-(\tfrac{t^2}2 a_1 a'_1) 
  c'(\th)+(ta_1(\th)+\frac{t^2}{2} a_2(\th))ic'(\th)  \mod t^3
\end{align*}
then $\langle c', c_t \rangle \equiv 0 \mod t^2$, hence this 2-jet is horizontal
and $\langle c_{tt}, ic' \rangle \equiv a_2 \mod t$ as required.

\subsection*{\nmb.{4.6}. Computation of the sectional curvature in 
$B_{i,f}(S^1_L,\mathbb R^2)$ at $C$}
We now go further. We use the following formula which is valid in a chart:
\begin{align}
&2R_f(m,h,m,h)=2G^A_f(R_f(m,h)m,h) =
\AMStag{\nmb:{1}}\\&
= -2d^2G^A(f)(m,h)(h,m)  +d^2G^A(f)(m,m)(h,h)  +d^2G^A(f)(h,h)(m,m)  
\notag\\&\quad
-2G^A(\Ga(h,m),\Ga(m,h))
+2G^A(\Ga(m,m),\Ga(h,h)) \notag
\end{align}
The sectional curvature at the two-dimensional subspace $P_f(m,h)$ of the
tangent space which is spanned by $m$ and $h$ is then given by:
\begin{equation}
k_f(P(m,h)) = - \frac{G^A_f(R(m,h)m,h)}{\|m\|^2\|h\|^2-G^A_f(m,h)^2}.
\AMStag{\nmb:{2}}
\end{equation}
We compute this directly for $f=0$. 
 From the expansion up to order 2 of $G^A_f(h,k)$ in \nmb!{4.5.1} we get:
\begin{multline}
\frac1{2!}d^2G^A(0)(m,l)(h,k) 
=\int_{S^1_L}hk\,\biggl(-\tfrac12m'l'+
\\
+A\Bigl(2(ml''+m''l)\ka^2
+ml\ka^4+\tfrac12m'l'\ka^2 +(ml'+m'l)\ka\ka' +m''l''\Bigr)
\biggr)\;d\th
\AMStag{\nmb:{3}}\end{multline}
Thus we have:
\allowdisplaybreaks\begin{align*}
&-d^2G^A(0)(m,h)(h,m)+\tfrac12 d^2G^A(0)(m,m)(h,h)+\tfrac12 d^2G^A(0)(h,h)(m,m) =
\\&=
%-2d^2G^A(0)(m,h)(h,m) = 
-2\int_{S^1_L}hm\,\Bigl(-\tfrac12m'h'+
\\&\qquad
+A\Bigl(2(mh''+m''h)\ka^2
+mh\ka^4+\tfrac12m'h'\ka^2 +(mh'+m'h)\ka\ka' +m''h''\Bigr)
\Bigr)\;d\th 
\\&\quad
%d^2G^A(0)(m,m)(h,h) = 
+\int_{S^1_L}hh\,\Bigl(-\tfrac12{m'}^2
+A\Bigl(4mm''\ka^2
+m^2\ka^4+\tfrac12{m'}^2\ka^2 +2mm'\ka\ka' +{m''}^2\Bigr)
\Bigr)\;d\th
\\&\quad
%d^2G^A(0)(h,h)(m,m) = 
+\int_{S^1_L}mm\,\Bigl(-\tfrac12h'h'
%+\\&\qquad\qquad\quad
+A\Bigl(4hh''\ka^2
+hh\ka^4+\tfrac12h'h'\ka^2 +2hh'\ka\ka' +h''h''\Bigr)
\Bigr)\;d\th
\\&=
\int_{S^1_L}\,\Bigl(
\tfrac12(A\ka^2-1)(mh'-m'h)^2 +A(mh''-m''h)^2 
\Bigr)\;d\th.
\end{align*}
For the second part of the curvature we have
\begin{align*}
&-G_0(\Ga_0(h,m),\Ga_0(m,h))+G_0(\Ga_0(m,m),\Ga_0(h,h)) =
\\&=
\int_{S^1_L}-\Bigl((\tfrac12\ka-\tfrac12A\ka^3 +A\ka'')hm 
  +2A\ka'(h'm+m'h) +2A\ka h'm'\Bigr)^2\frac{d\th}{1+A\ka^2}
\\&\quad
+\int_{S^1_L}\Bigl((\tfrac12\ka-\tfrac12A\ka^3 +A\ka'')m^2
  +4A\ka'm{m'} +2A\ka {m'}^2\Bigr)
\\&\qquad\qquad
\Bigl((\tfrac12\ka-\tfrac12A\ka^3 +A\ka'')h^2 
  +4A\ka'hh' +2A\ka {h'}^2\Bigr)
  \frac{d\th}{1+A\ka^2}
\\&=
\int_{S^1_L}\Bigl(
(A\ka^2-A^2\ka^4 +2A^2\ka\ka''-4A^2{\ka'}^2)(mh'-m'h)^2
\Bigr)\frac{d\th}{1+A\ka^2}
\end{align*}
Thus we get
\begin{align*}
&R_0(m,h,m,h)=G^A_0(R_0(m,h)m,h) =
\\&=
\int_{S^1_L}\,\Bigl(
\tfrac12(A\ka^2-1)(mh'-m'h)^2 +A(mh''-m''h)^2 
\Bigr)\;d\th
\\&\quad
+\int_{S^1_L}\Bigl(
(A\ka^2-A^2\ka^4 +2A^2\ka\ka''-4A^2{\ka'}^2)(mh'-m'h)^2
\Bigr)\frac{d\th}{1+A\ka^2}
\end{align*}
Letting $W = mh'-hm'$ be the Wronskian of $m$ and $h$ and simplifying, we have:
\begin{equation}
\boxed{\quad\begin{aligned}
R_0&(m,h,m,h)=
\\&
=\int_{S^1_L}\,\Bigl(\frac{
-(A\ka^2-1)^2+4A^2\ka\ka''-8A^2{\ka'}^2}{2(1+A\ka^2)}\Bigr)W^2
d\th
+\int_{S^1_L}A{W'}^2 d\th
\end{aligned}\quad}
\AMStag{\nmb:{4}}
\end{equation}
What does this formula say? First of all, if
$\text{supp}(m) \cap \text{supp}(h) = \emptyset$, the sectional curvature in
the plane spanned by $m$ and $h$ is 0. Secondly, we can divide the curve $c$
into two parts:
\begin{align*}
c^+_A &= \text{ set of points where } \ka \ka'' < 2(\ka')^2 +
\bigl(\tfrac{A^{-1}-\ka^2}2 \bigr)^2 \\
c^-_A &= \text{ set of points where } \ka \ka'' > 2(\ka')^2 +
\bigl(\tfrac{A^{-1}-\ka^2}2 \bigr)^2.
\end{align*}
Note that if $A$ is sufficiently small, $c^-_A = \emptyset$ and even if $A$
is large, $c^-_A$ need not be non-empty. 
But if $\text{supp}(m),\text{supp}(h) \subset c^-_A$,
the sectional curvature is always negative.
The interesting case is when $\text{supp}(m),\text{supp}(h) \subset c^+_A$. 
We may introduce the
self-adjoint differential operator on $L^2(S^1)$:
$$Sf = f'' +\frac{(A\ka^2-1)^2-4A^2\ka\ka''+8A^2{\ka'}^2}{2A(1+A\ka^2)}f$$
so that $R = -A\langle SW,W \rangle$.
The eigenvalues of $S$ tend to $-\infty$, hence $S$ has a finite number of
positive eigenvalues. If we take, for example, $m=1$ and $h$ such that $h'$
is in the span of the positive eigenvalues, the corresponding sectional
curvature will be positive. 
In general, the condition that the sectional curvature be
positive is that the Wronskian $W$ have a sufficiently large component in
the positive eigenspace of $S$. The special case where $c$ is the unit circle
may clarify the picture: then
$$ Sf = f'' + \frac{(A-1)^2}{2A(1+A)}f$$
and the eigenfunctions are linear combinations of sine's and cosine's. 
It is easy to see that for any
$A$, a plane spanned by $m$ and $h$ of pure frequencies $k$ and $l$ will have
positive curvature if and only if $k$ and $l$ are sufficiently near each
other (asymptotically $|k-l| < |A-1|/\sqrt{A+a^2}$), hence `beat' at a low
frequency.

\subsection*{\nmb.{4.7}. The sectional curvature for the induced $H^0$-metric 
on $B_{i,f}(S^1_L,\mathbb R^2)$ in a chart }
In the setting of \nmb!{4.2} we have for 
$f\in C^\infty(S^1_L,(-\ep,\ep))$ and $h,k\in C^\infty(S^1_L,\mathbb R)$
\allowdisplaybreaks\begin{align}
G^0_f(h,k) &= ((\pi\o\ps)^*G^0)_f(h,k) = 
G^0_{\pi(\ps(f))}\Bigl(T_f(\pi\o\ps)h,T_f(\pi\o\ps)k\Bigr)
\AMStag{\nmb:{1}}\\
&=G^0_{\ps(f)}\Bigl((T_f\ps.h)^\bot,(T_f\ps.k)^\bot\Bigr)
\notag\\&
=\int_{S^1_L}\frac{hk(1-f\ka_c)^2}{\sqrt{(1-f\ka_c)^2+{f'}^2}}\;d\th
\notag\end{align}
At the center of the chart described in \nmb!{4.4}, i.e., for $f=0$, the
Christoffel symbol \nmb!{4.5.2} for $A=0$ becomes
\begin{equation}
\Ga_0(h,k) = \tfrac12\ka_c hk 
\AMStag{\nmb:{2}}
\end{equation}
The curvature \nmb!{4.6.4} at $f=0$ for $A=0$ becomes 
\begin{align}
R_0(m,h,m,h)&=G_0(R_0(m,h)m,h) =
\notag\\&
=-\tfrac12\int_{S^1_L} (h'm-hm')^2\,d\th 
=-\tfrac12\int_{S^1_L} W(m,h)^2\,d\th 
\AMStag{\nmb:{3}}
\end{align}
and the sectional curvature $k_0(P(m,h))$ from \nmb!{4.5.2} for $A=0$ and
$f=0$
%\begin{displaymath}
%k(P(m,h)) = - \frac{G^0_0(R(m,h)m,h)}{\|m\|^2\|h\|^2-G_f(m,h)^2}.
%\end{displaymath}
is non-negative.

In the full chart \nmb!{4.2}, starting from the metric \nmb!{4.6.1},
we managed to compute the full geodesic equation not just for $f=0$ but 
for general $f$, so long as $A=0$. The outcome is
\begin{align}
\Ga_f(h,h)
&= \frac{\ka_c h^2}{1-f\ka_c}
  +\frac{-\tfrac12\ka_c(1-f\ka_c) h^2 + (\tfrac12 h^2f''+2hh'f')}
  {\bigl((1-f\ka_c)^2+{f'}^2\bigr)}
\notag\\&\qquad
-\frac{\ka_c h^2{f'}^2}
  {(1-f\ka_c)\bigl((1-f\ka_c)^2+{f'}^2\bigr)}
+ \frac{\tfrac32\ka_c(1-f\ka_c)h^2{f'}^2-\tfrac32 h^2{f'}^2f''}
  {\bigl((1-f\ka_c)^2+{f'}^2\bigr)^{2}}.
\AMStag{\nmb:{4}}
\end{align}
The geodesic equation is thus
\begin{align}
f_{tt}
&= -\frac{\ka_c f_t^2}{1-f\ka_c}
  -\frac{-\tfrac12\ka_c(1-f\ka_c) f_t^2 
          + (\tfrac12 f_t^2f_{\th\th}+2f_tf_{t\th}f_\th)}
  {\bigl((1-f\ka_c)^2+{f_\th}^2\bigr)}
\notag\\&\qquad
+\frac{\ka_c f_t^2{f_th}^2}
  {(1-f\ka_c)\bigl((1-f\ka_c)^2+{f_\th}^2\bigr)}
- \frac{\tfrac32\ka_c(1-f\ka_c)f_t^2{f_\th}^2-\tfrac32f_t^2{f_\th}^2f_{\th\th}}
  {\bigl((1-f\ka_c)^2+{f_\th}^2\bigr)^{2}}.
\AMStag{\nmb:{5}}
\end{align}
For $A>0$ we were unable to get the analogous result.

\section*{\totoc\nmb0{5}. 
Examples and numerical results}

\subsection*{\nmb.{5.1}. 
The geodesics running through concentric circles}
The simplest possible geodesic in $B_i$ is given by the set of all 
circles with common center. Let $C_r$ be the circle of radius $r$ with
center the origin. Consider the path of such circles $C_{r(t)}$ given by the 
parametrization $c(t,\th)=r(t)e^{i\th}$, where $r(t)$ is a smooth increasing 
function $r:[0,1]\to \mathbb R_{>0}$. Then $\ka_c(t,\th)=1/r(t)$. If we vary 
$r$ then the horizontal energy and the variation of this curve are
\begin{align*}
E^{\text{hor}}_{G^A}(c) &= \frac12\int_0^1\int_{S^1}
\bigl(1+A/r^2\bigr)r_t^2 r\,d\th\,dt
\\
\p_s|_{s=0}E^{\text{hor}}_{G^A}(c) 
  &= \int_0^1\int_{S^1}\Bigl(1+\frac{A}{r^2}\Bigr)
  r_s\Bigl(-r_{tt}-\frac{(1-A/r^2)}{2(r+A/r)}r_t^2\Bigr)
  r\,d\th\,dt
\end{align*}
so that $c$ is a geodesic if and only if
\begin{equation}
r_{tt}+\frac{(1-A/r^2)}{2(r+ A/r)}r_t^2 = 0.
\AMStag{\nmb:{1}}
\end{equation}
Also the geodesic equation \nmb!{4.1.1} reduces to \thetag{\nmb|{1}} for
$c$ of this form. 

The solution of \thetag{\nmb|{1}} can be written in terms
of the inverse of a complete elliptic integral of the second kind.
More important is to look at what happens for small and large $r$. As
$r \rightarrow 0$, the ODE reduces to:
$$ r_{tt} -\frac{r_t^2}{2r} = 0$$
whose general solution is $r(t) =C(t-t_0)^2$ for some contants $C,t_0$.
In other words, at one end, the path ends in finite time with the circles
imploding at their common center. Note that $r' \rightarrow 0$ as 
$r \rightarrow 0$ but not fast enough to prevent the collapse. 
On the other hand, as $r \rightarrow \infty$, the ODE becomes:
$$ r_{tt} +\frac{r_t^2}{2r} = 0$$
whose general solution is $r(t) = C (t-t_0)^{2/3}$ for some constants
$C,t_0$. Thus at the other end of the geodesic, the circles expand forever
but with decreasing speed.

An interesting point is that this geodesic has conjugate points on it, so that
it is a extremal path but not a local minimum for length over all intervals. 
This is a concrete reflection
of the collapse of the metric when $A=0$. To work this out, take any $f(\th)$
such that $\int_0^{2\pi} f d\th = 0$ and any function $a(t)$. Then
$X = f(\th)a(t)\partial/\partial r$ is a vector field along the geodesic, i.e.\
a family of tangent vectors to $B_e$ at each circle $C_{r(t)}$ normal to
the tangent to the geodesic. Its length is easily seen to be:
$$ 
\|X\|^2_{C_{r(t)}} = \Bigl(r(t)+\frac{A}{r(t)}\Bigr)
a(t)^2\int_0^{2\pi} f(\th)^2 d\th.
$$
We need to work out its covariant derivative:
$$ \nabla_\frac{d}{dt} (X) = f(\th) a_t \frac{\partial}{\partial r} +
\Gamma_{C_r} \bigl(r_t \frac{\partial}{\partial r}, 
f(\th) a \frac{\partial}{\partial r}\bigr).$$
Using a formula for the Christoffel symbol which we get from \nmb!{4.2.2}
by polarizing, and noting that
$\ka \equiv 1/r, \ka_s \equiv 0$, we get:
\begin{align*}
\nabla_\frac{d}{dt} (X) &= f(\th) a_t \frac{\partial}{\partial r} +
f(\th) a r_t \Bigl(\frac{1-A/r^2}{2(r+A/r)}\Bigr) \frac{\partial}{\partial r} \\
&= f(\th) \bigl(r+A/r\bigr)^{-1/2} 
  \Bigl( \bigl(r+A/r\bigr)^{1/2} a\Bigr)_t \frac{\partial}{\partial r}.
\end{align*}
(This formula also follows from the fact that the vectors 
$(r+A/r)^{-1/2} \partial/\partial r$ have length independent of $t$, hence
covariant derivative zero.) Jacobi's
equation is therefore:
\begin{equation} 
\AMStag{\nmb:{2}}
f(\th) \bigl( r+A/r \bigr)^{-1/2} \Bigl( \bigl(r+A/r\bigr)^{1/2} 
a\Bigr)_{tt} \frac{\partial}{\partial r} +
R\bigl(X,r_t \frac{\partial}{\partial r}\bigr)
\bigl(r_t \frac{\partial}{\partial r}\bigr) = 0,
\end{equation}
where $R$ is the curvature tensor. For later purposes, it is convenient to
write this eqation using $r$ as the independent variable along the geodesic 
rather than $t$ and think of $a$ as a function of $r$. Note that for any function
$b$ along the geodesic, $b_t = b_r r_t$ and
$$b_{tt} = b_{rr} r_t^2 + b_r r_{tt} = 
\bigl(b_{rr}-\frac{(1-A/r^2)}{2(r+A/r)} b_r \bigr) r_t^2.$$
Then a somewhat lengthy bit of algebra shows that:
\begin{align*}
\bigl( r+A/r \bigr)^{-\tfrac12} \Bigl( \bigl(r+A/r\bigr)^{\tfrac12} 
a\Bigr)_{tt}  &=
\bigl( r+A/r \bigr)^{-\tfrac14} \Bigl( \bigl(r+A/r\bigr)^{\tfrac14} 
a\Bigr)_{rr}r_t^2 + F(r) a r_t^2, \\ \text{where} \qquad
F(r) &= -\frac{5}{16} \Bigl( \frac{1-A/r^2}{r+A/r} \Bigr)^2 +
\frac{A}{2r^3(r+A/r)}.
\end{align*}

To work out the structure of $R$ in this case, 
use the fact that the circles $C_r$ and the vector field 
$\partial/\partial r$
are invariant under rotations. This means that the map $f \mapsto
R(\partial/\partial r, f \partial/\partial r)(\partial/\partial r)$
has the two properties: it commutes with rotations and it is symmetric. The
only such maps are diagonal in the Fourier basis, i.e.\ there are real
constants $\lambda_n$ such that 
$$R\Bigl(\partial/\partial r, \Bigl\{ \begin{matrix}
\cos(n\th)\partial/\partial r\\ \sin(n\th)\partial/\partial r \end{matrix}\Bigr)
(\partial/\partial r) = \lambda_n \Bigl\{ 
\begin{matrix} \cos(n\th)\partial/\partial r \\ \sin(n\th)\partial/\partial r 
\end{matrix}.$$
To evaluate $\lambda_n$, we take the inner product with $\cos(n\th)$ (or 
$\sin(n\th)$ and use our calculation of $R_0(m,h,m,h)$ in section
\nmb!{4.6} to show:
\begin{align*}
\Bigl\langle R\bigl(\frac{\partial}{\partial r}, 
\cos(n\th)\frac{\partial}{\partial r}\bigr)
(\frac{\partial}{\partial r}), \cos(n\th) \frac{\partial}{\partial r} 
\Bigr\rangle 
&=
R_0\Bigl(\frac{\partial}{\partial r}, \cos(n\th) \frac{\partial}{\partial r}, 
\frac{\partial}{\partial r}, \cos(n\th) \frac{\partial}{\partial r}\Bigr) \\
&= \int_0^{2\pi}\Big( -\frac{(1-A/r^2)^2}{2(1+A/r^2)} W^2 + A W'^2 \Big)r d\th
\end{align*}
where 
$$
W = 1.\frac{d}{ds}\cos(n\th) = -n\frac{\sin(n\th)}{r}
\quad\text{  and }\quad
W'=\frac{d}{ds}W = -n^2\frac{\cos(n\th)}{r^2}.
$$
Simplifying, this gives:
\begin{align*}
\lambda_n \|\cos(n\th)\frac{\partial}{\partial r}\|^2 &=
\int_0^{2\pi}\Bigl(-\frac{(1-A/r^2)^2}{2(r+A/r)} n^2 \sin^2(n\th)
+ \frac{A}{r^3} n^4 \cos^2(n\th) \Bigr) d\th \\
&=  -\frac{(1-A/r^2)^2}{2(r+A/r)} n^2 \pi + \frac{A}{r^3} n^4 \pi
\end{align*}
hence
$$ \lambda_n = -\frac{(1-A/r^2)^2}{2(r+A/r)^2} \cdot n^2 +
\frac{A}{r^3(r+A/r)} \cdot n^4.$$
Thus for $X = \cos(n\th) a_n(t) \partial/\partial r$, if we combine
everything, Jacobi's equation reads:
\begin{equation}
\AMStag{\nmb:{3}}
\begin{aligned}
 \bigl( r+A/r \bigr)^{-\tfrac14}& \Bigl( \bigl(r+A/r\bigr)^{\tfrac14} 
a_n \Bigr)_{rr} =
\\&
= \Bigl(-\frac{(1-A/r^2)^2}{2(r+A/r)^2} (n^2-\tfrac58) +
\frac{A}{r^3(r+A/r)} (n^4 - \tfrac12)\Bigr) a_n.
\end{aligned}
\end{equation}

\begin{figure}[htb]
\begin{center}
\epsfig{width=4in,file=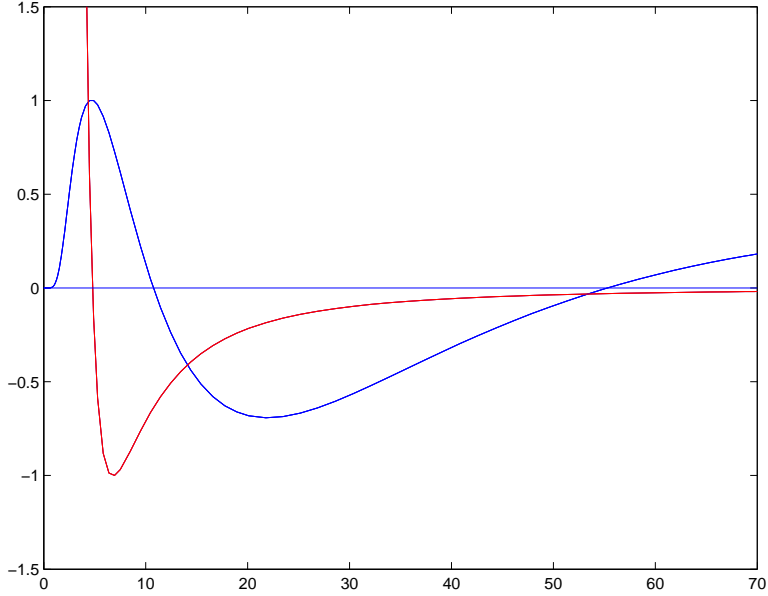}
\caption{The potential in the Jacobi ODE and its solution for an 
infinitesimal triangular perturbation of the circles in the geodesic of
concentric circles. Note the first conjugate point at $10.77 \sqrt{A}$.}
\end{center}
\end{figure}

Calling the right hand side the {\it potential} of Jacobi's equation, 
we can check that for each $n$, the potential is positive for small $r$, 
negative for large $r$ and it has one zero, approximately at $\sqrt{2A} n$ 
for large $n$. Thus, for small $r$, these perturbations diverge from the 
geodesic of circles. For large $r$, if we write 
$b_n = \bigl(r+A/r\bigr)^{1/4} a_n$, then
Jacobi's equation approaches:
$$ (b_n)_{rr} \approx -\frac{n^2 - 0.625}{2r^2} b_n.$$
This is solved by $b_n = c x^\lambda + c' x^{\lambda'}$ where $\lambda, \lambda'$ 
are solutions of $\lambda^2 -\lambda = -(n^2-0.625)/2$. For $n=1$, 
$\lambda, \lambda'$ are real and $b_n$ has no zeros; but for $n>1$,
$\lambda, \lambda'$ have an imaginary part, say $i\gamma_n$, and
$$
b_n \approx \sqrt{r} \big(c \cos(\gamma_n\log(r)) 
  + c' \sin(\gamma \log(r))\big)
$$
with infinitely many zeros.

Figure 4 shows the solution for $n=3$ which approaches 0 as $r \rightarrow 0$.
The first zero of this solution is about $10.77 \sqrt{A}$, making it a conjugate
point of $r=0$. For other $n$, the first such conjugate point appears to be
bigger, so we conclude: on any segment $0 < r_1 < r_2 < 10.77 \sqrt{A}$, the
geodesic of circles is locally (and presumably globally) minimizing. 
              
\subsection*{\nmb.{5.2}. The geodesic connecting two distant curves}
For any two distant curves $C_1,C_2$, one can construct paths from one 
to the other by (a) changing $C_1$ to some auxiliary curve $D$ near $C_1$, 
(b) translating $D$ without modifying it to a point near $C_2$ and 
(c) changing the translated curve $D$ to $C_2$. 
If $C_1$ and $C_2$ are very far from each other, 
the energy of the translation will dominate the energy required to modify 
them both to $D$. Thus we expect that a geodesic between distant curves will 
aymptotically utilize a curve $D$ which is optimized for least energy translation. 
To find such curves $D$, heuristically we may argue that it should be a curve 
such that the path given by all its translates in a fixed direction is a geodesic.

Such geodesics can be found as special cases of the general geodesic. We fix $e = 
(1,0)$ as the direction of translation and assume that the path $\{D+te\}$ is a 
geodesic. We need to express this geodesic up to order $O(t^2)$ in the chart used 
in section \nmb!{4.4}. Let $c(s)$ be arc length parametrization of $D$ and $\th(s)$ be 
the orientation of $D$ at point $c(s)$, i.e.\ $c_s = \cos(\th) + i\sin(\th)$.
Then a little calculation shows that if we reparametrize nearby curves
via $\tilde s = s - \langle e, c_s \rangle t$, then the path of translates in 
direction $e$ is just:
\begin{align*} 
c(\tilde s) + t e &= c(s) + \bigl( t\langle e,ic_s \rangle + 
\frac{t^2}{2} \langle e, c_s \rangle^2 \ka + O(t^3)\bigr) ic_s \\
&= c(s) +\bigl( - \sin(\th(s)) t +\frac{t^2}{2} \cos^2(\th(s) \ka) + 
O(t^3)\bigr) ic_s.
\end{align*}
Thus, in the notation of \nmb!{4.2}, $a=-\sin(\th)$, hence
$a_s = -\cos(\th) \ka$ and, moreover, $a_t = 
\cos^2(\th) \ka$. Substituting this in the geodesic formula \nmb!{4.2.1}, 
we get
\begin{multline*} (1+A\ka^2) \cos^2(\th)\ka = 
\\
=\frac{\ka \sin^2(\th)}{2} + A\bigl(
(\ka_{ss} - \frac{\ka^3}{2})\sin^2(\th) + 4\cos(\th)\sin(\th)\ka \ka_s
+ 2 \ka^3 \cos^2(\th) \bigr).
\end{multline*}
Since $\ka = \th_s$, this becomes, after some manipulation, a singular 
third order equation for $\th(s)$:
$$ \th_{sss} = 4\cot(\th) \th_s \th_{ss} +
(\tfrac12 - \cot^2(\th))\th_s (\th_s^2 - \tfrac1{A}).$$
One solution of this equation is $\th(s) \equiv \tfrac1{\sqrt{A}}$, i.e.\ a 
circle of radius $\sqrt{A}$. In fact, this seems to be the only simple closed
curve which solves this equation. However, if we drop smoothness, a weak
solution of this equation is given by the $C^1$, piecewise $C^2$-curve made 
up of 2 semi-circles of radius $\sqrt{A}$ joined by 2 straight line segments 
parallel to the vector $e$ and separated by the distance $2\sqrt{A}$ (as in 
figure 5). Note that such `cigar'-shaped curves can be made with line segments
of any length.

\begin{figure}[htb] \begin{center}
\epsfig{width=4in,file=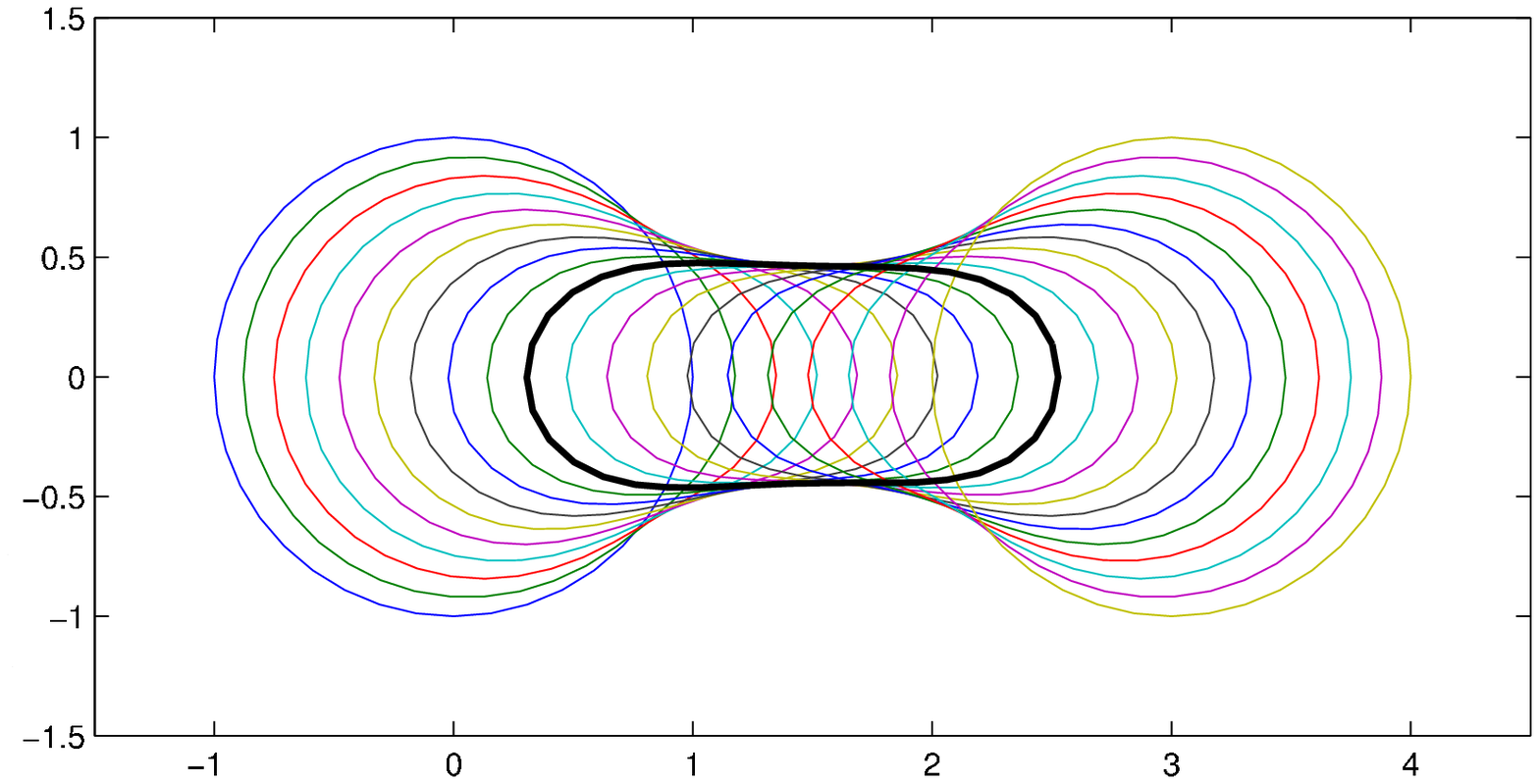}\\
\epsfig{width=4in,file=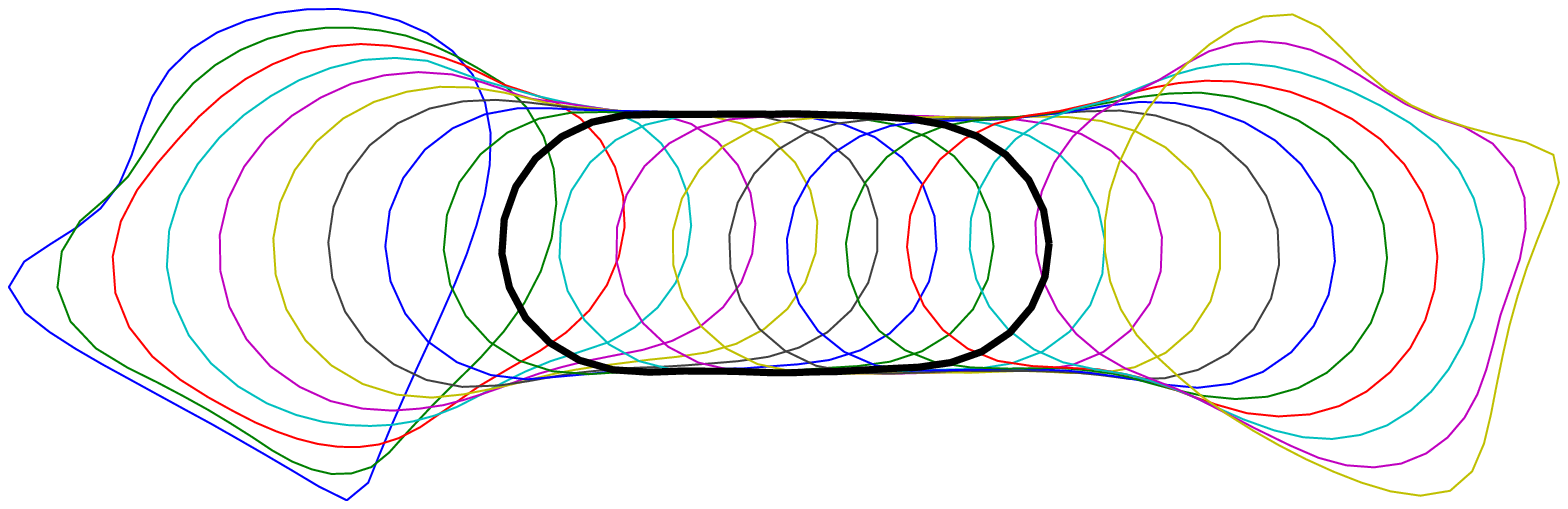}
\caption{On the top, the geodesic joining circles of radius 1 at distance 3 apart
with $A=.1$ (using 20 time samples and a 40-gon for the circle). On the
bottom, the geodesic joining 2 `random' shapes of size about 1 at distance
5 apart with $A=.25$ (using 20 time samples and a 48-gon approximation for
all curves). In both cases the middle curve which is highlighted.}
\end{center} \end{figure}

A numerical approach to minimize $E^{\text{hor}}_{G^1}(c)$ for variations
$c$ with initial and end curves circles at a certain distance produced the
two such geodesics shown in Figure 5. Note that the middle curve is indeed
close to such a `cigar'-shape. However, the width of this shape is somewhat
greater than $2\sqrt{A}$: this is presumably because the endcurves of this
path are not sufficiently far apart. Thus experiments as well as the theory 
suggest strongly that geodesics joining any two curves sufficiently far apart 
compared to their size asymptotically approach a constant `cigar'-shaped 
$C^1$-intermediate curve made up of 2 semi-circles of radius $\sqrt{A}$ and 
2 parallel line segments. We conjecture that this is true.

\subsection*{\nmb.{5.3} The growth of a `bump' on a straight line, when $A=0$}
We have seen above that the geodesic spray is locally well-defined when $A=0$.
To understand this spray and see whether it appears to have global solutions,
we take that the initial curve contains a segment with curvature identically
zero, i.e.\ contains a line segment, and that the initial velocity $a$ is set
to a smooth function with compact support contained in this segment. For
simplicity, we take the velocity $a$ to be a cubic B-spline, i.e.\ a piecewise
cubic which is $C^2$ with 5 non-$C^3$ knots approximating a Gaussian blip. 
The result of integrating is shown in Figure 6. Note several things: first,
where the curvature is zero, the curve moves with constant velocity if we
follow the orthogonal trajectories. Secondly, where the curve is moving
opposite to its curvature (like an expanding circle, the part in the middle), 
it is deccelerating; but where it is moving with its curvature (like a contracting
circle, the parts on the 2 ends), it is accelerating. This acceleration in the
2 ends, creates higher and higher curvature until a corner forms. In the
figure, the simulation is stopped just before the curvature explodes. In the
middle, the curve appears to be getting more and more circular. As the corners
form, the curve is approaching the boundary of our space. Perhaps, with the right
entropy condition, one can prolong the solution past the corners with a suitable
piecewise $C^1$-curve.

\begin{figure} \begin{center}
\epsfig{width=3in,file=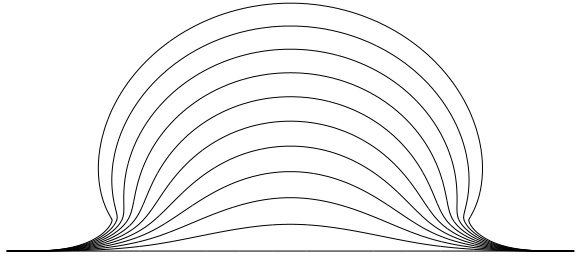}
\caption{The forward integration of the geodesic equation when $A=0$, 
starting from a straight line in the direction given by a smooth bump-like 
vector field. Note that two corner like singularities with curvature going
to $\infty$ are about to form.}
\end{center} \end{figure}

Although this calculation assumes $A=0$, one will find very similar geodesics
when $A$ is much smaller than $1/\ka^2, 1/(\ka_s \log(a)_s)$ and $\ka/\ka_{ss}$,
so that the dominant terms in the geodesic equation are those without an $A$.
In other words, geodesics between large smooth curves are basically the
same as those with $A=0$.

\subsection*{\nmb.{5.4} Several geodesic triangles in $B_e$}
We have examined dilations, translations and the evolution of blips. We look
next at rotations. To get a pure rotational situation, we consider ellipses
centered at $(0,0)$ with the same eccentricity 3 and maximum radius 1, but 
differently oriented. We take 3 such ellipses, with orientations differing
by $60^\circ$ and $120^\circ$ degrees. Joining each pair by a geodesic, we get
a triangle in $B_e$. 

We wanted to examine whether along the geodesic joining 2 such ellipses
(a) one ellipse rotates into the other or (b) the initial ellipse shrinks
towards a circle, while the final ellipse grows, independently of one
another. It turns out that, depending on the value of $A$, both can
happen. Note that we get similar geodesics by either changing $A$ or making
the ellipses smaller or larger with $A$ held fixed. For each $A$, we get
an absolute distance scale with unit $1/\sqrt{A}$ and, if the ellipses are
bigger than this, (b) dominates, while, if smaller, (a) dominates.

The results are shown in Figure 7. 
We have taken the three values $A=1,0.1$ and $0.01$. 
For each value, on the top, we show the geodesic joining 2 of the
ellipses as a sequence of curves in their common ambient $\mathbb R^2$. Below this, 
we show the triple of geodesics as a triangle, by displaying the intermediate
curves as small shapes along lines joining the ellipses. This Euclidean triangle
is being used purely for display, to indicate that the computed structure is
a triangle in $B_e$. Note that for $A=1$, the intermediate shapes are very close 
to ellipses, whose axes are rotating; 
while for $A=0.01$, the bulges in one ellipse shrink while those of the other grow. 

\begin{figure} \begin{center}
\epsfig{width=\columnwidth,file=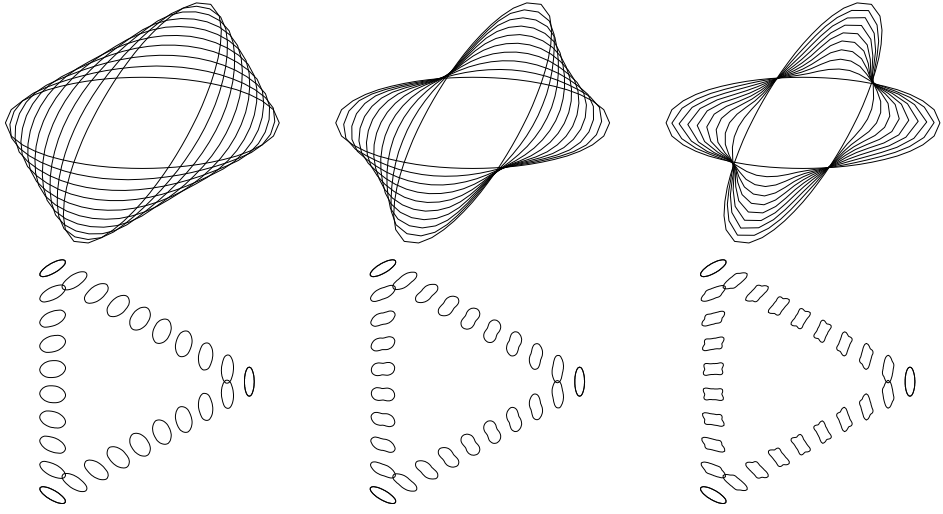} 
\caption{Top Row: Geodesics in three metrics joining the same two ellipses.
The ellipses have eccentricity 3, the same center and are
at $60^\circ$ degree angles to each other. At left, $A=1$;
in middle $A=0.1$; on right $A=0.01$. Bottom Row: Geodesic triangles
in $B_e$ formed by joining three ellipses at angles 0, 60 and 120 degrees,
for the same three values of $A$. Here the intermediate shapes are just
rotated versions of the geodesic in the top row but are laid
out on a plane triangle for visualization purposes.}
\end{center} \end{figure}

We can also compute the angles in $B_e$ between the sides of this triangle.
They work out to be $34^\circ$ when $A=1$, i.e.\ the angle sum for the triangle 
is $102^\circ$, much less than $\pi$ radians, showing strong negative sectional 
curvature in the plane containing this triangle. But if $A=0.1$ or $0.01$, 
the angle is $77^\circ$ and $69^\circ$ respectively, giving more than
$\pi$ radians in the triangle. Thus the sectional curvature is positive 
for such small values of $A$.

\subsection*{\nmb.{5.5} Notes on the numerical simulations}
All simulations in this paper were carried out in MatLab. The forward
integration for the geodesic equation for $A=0$ was carried out by the
simplest possible finite difference scheme. This seems very stable and
reliable. Solving for the geodesics was done using the MatLab minimization
routine {\tt fminunc} using both its medium and large scale modes. 
This, however, was quite unstable due to discretization
artifacts. A general path between two curves was represented by a matrix
of points in $\mathbb R^2$, approximating each curve by a polygon and sampling
the path discretely. The difficulty is that when the polygons have very acute
angles, the discretization tends to be highly inaccurate because of the
high curvature localized at one vertex. Initially, in order to minimize
the number of variables in the problem, we tried to use small numbers of
samples and higher order accurate discrete approximations to the derivatives. 
In all these attempts, the discrete approximation ``cheated'' by finding 
minima to the energy of the path with polygons with very small angles. The
only way we got around this was to use first order accurate expressions for
the derivatives and relatively large numbers of samples (e.g.\ 48 points
on each curve, 20 samples along the geodesic, hence 
$2\times 20 \times 48 = 1920$ variables in the expression for the energy.

Another problem is that the energy only depends on the path of unparametrized
curves and is independent of the parametrization. To solve this, we added a
term to the energy which is minimized by constant speed parametrizations.
This still leaves a possibly wandering basepoint, and we added $\epsilon$
times another term which asked that all points on each curve should move
as normally as possible. In practice, if the initialization was reasonable,
this term was not needed. The final discrete energy that was minimized was
this. Let $x_{i,j}$ be the $i^{th}$ sample point on the $j^{th}$ curve $C_j$.
For each $(i,j)$, estimate the sum of the squared curvature of $C_j$
plus the squared acceleration of the parametrization by:
$$ k(i,j) = \tfrac12 \Bigl( \frac{1}{\|x_{i-1,j}-x_{i,j}\|^4} + 
\frac{1}{\|x_{i,j}-x_{i+1,j}\|^4} \Bigr) \cdot
\|x_{i-1,j} - 2 x_{i,j} + x_{i+1,j}\|^2.$$
(The harmonic mean of the segment lengths is used here to further force
the parametrization to be uniform.) Then, for each $(i,j)$, the {\it four}
triangles $t = \{a=(i,j), b=(i\pm 1,j), c=(i,j\pm 1)\}$ around $(i,j)$ are
considered and the energy is taken to be:
$$ \sum_{i,j,t} \Bigl( \frac{\langle (x_a - x_b), (x_a - x_c)^\perp \rangle^2
+\epsilon \langle (x_a-x_b),(x_a-x_c) \rangle^2}{\|x_a-x_b\|}\Bigr)
(1+A k(a)).$$
We make no guarantees about the accuracy of this simulation! The results, 
however, seem to be stable and reasonable.
%\newpage

%\nocite{*}\bibliographystyle{plain}\bibliography{\jobname}

\end{document}